\documentclass{article}

\usepackage{arxiv}

\usepackage[utf8]{inputenc} 
\usepackage[T1]{fontenc}    
\usepackage{hyperref}       
\usepackage{url}            
\usepackage{booktabs}       
\usepackage{amsfonts}       
\usepackage{nicefrac}       
\usepackage{microtype}      
\usepackage{lipsum}		
\usepackage{graphicx}
\usepackage{natbib}
\usepackage{doi}

\usepackage{amsmath,amssymb}
\usepackage{amsthm}
\usepackage{bm,upgreek}
\usepackage{epstopdf}
\usepackage[FIGBOTCAP,TABTOPCAP,bf,tight]{subfigure}
\usepackage{setspace}
\usepackage{enumerate,etaremune}
\usepackage{tabularx,multirow}
\usepackage{framed}
\usepackage{hhline}
\usepackage{cancel}
\usepackage{empheq}
\usepackage{tikz}
\usepackage[title]{appendix}
\usetikzlibrary{shapes.geometric}
\usetikzlibrary{shapes,arrows}
\usetikzlibrary{plotmarks}
\usetikzlibrary{calc}

\newcommand{\tensor}[1]{\bm{#1}}

\newcommand{\transpose}[1]{{#1}^{\mathsf T}}

\newcommand{\rn}[1]{\uppercase\expandafter{\romannumeral #1\relax}}
\DeclareMathOperator{\grad}{\nabla}

\DeclareMathOperator{\diver}{\nabla\cdot}

\DeclareMathOperator{\tr}{tr}

\theoremstyle{remark}

\renewcommand{\vec}[1]{\ensuremath{\boldsymbol{#1}}}

\renewcommand{\transpose}[1]{{#1}^{\mathsf T}}

\title{Can ChatGPT implement finite element models for geotechnical engineering applications?}

\date{} 

\author{ 
    Taegu Kim\\
	School of Civil and Environmental Engineering\\
	Yonsei University\\
	Seoul 03722, Republic of Korea \\
	\texttt{xorn6487@yonsei.ac.kr} \\
	\And
	\href{https://orcid.org/0000-0001-8701-543X}{\includegraphics[scale=0.06]{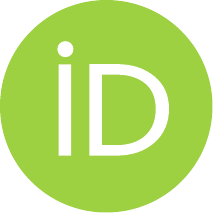}\hspace{1mm}Tae Sup Yun} \\
	School of Civil and Environmental Engineering\\
	Yonsei University\\
	Seoul 03722, Republic of Korea \\
	\texttt{taesup@yonsei.ac.kr} \\
	\And
	\href{https://orcid.org/0000-0002-5032-5255}{\includegraphics[scale=0.06]{orcid.pdf}\hspace{1mm}Hyoung Suk Suh} \\
	Department of Civil and Environmental Engineering\\
	Case Western Reserve University\\
	Cleveland, OH 44106, USA \\
	\texttt{hssuh@case.edu} \\
}

\begin{document}
\maketitle

\begin{abstract}
This study assesses the capability of ChatGPT to generate finite element code for geotechnical engineering applications from a set of prompts. 
We tested three different initial boundary value problems using a hydro-mechanically coupled formulation for unsaturated soils, including the dissipation of excess pore water pressure through fluid mass diffusion in one-dimensional space, time-dependent differential settlement of a strip footing, and gravity-driven seepage. 
For each case, initial prompting involved providing ChatGPT with necessary information for finite element implementation, such as balance and constitutive equations, problem geometry, initial and boundary conditions, material properties, and spatiotemporal discretization and solution strategies. 
Any errors and unexpected results were further addressed through prompt augmentation processes until the ChatGPT-generated finite element code passed the verification/validation test. 
Our results demonstrate that ChatGPT required minimal code revisions when using the FEniCS finite element library, owing to its high-level interfaces that enable efficient programming. 
In contrast, the MATLAB code generated by ChatGPT necessitated extensive prompt augmentations and/or direct human intervention, as it involves a significant amount of low-level programming required for finite element analysis, such as constructing shape functions or assembling global matrices. 
Given that prompt engineering for this task requires an understanding of the mathematical formulation and numerical techniques, this study suggests that while a large language model may not yet replace human programmers, it can greatly assist in the implementation of numerical models. 
\end{abstract}

\keywords{Finite element method \and Large language model \and Code generation \and Prompt engineering \and Unsaturated soil}

\section{Introduction}
\label{sec:intro}
In a wide range of geotechnical engineering applications from traditional tasks such as slope stability analysis or foundation design to modern challenges like CO$_2$ geological sequestration or nuclear waste disposal, accurate modeling of hydro-mechanical processes in unsaturated media is essential for reliable forward predictions \citep{matsui1992finite, yazdchi1999dynamic, hudson2005guidance, gens2010modelling, lei2015tough2biot}. 
This requires a strong coupling between the solid skeleton and fluids within the pores, as the presence of pore fluids may induce excess pore pressures that introduce rate-dependence to the overall mechanical responses, while deformation of the solid matrix can alter the geometry and topology of the pore network, significantly influencing the fluid pressures and saturation ratios therein. 
At the continuum scale, this hydro-mechanical coupling can be understood within the context of mixture theory, where each phase is considered as a constituent of a continuum mixture \citep{rice1976some, coussy1998mixture}. 
While the mixture theory provides a theoretical foundation for developing balance equations for multi-phase bodies from first principles, the resulting partial differential equations are typically not solvable using classical analytical methods. 
This challenge has led to the development of various techniques designed to solve coupled problems numerically \citep{de1996highlights, zienkiewicz1999computational, jeremic2008numerical, chen2023advanced}.

Although technically feasible, capturing the intricate interplay among different phase constituents is, nevertheless, not a trivial task. 
One possible strategy is to employ an operator splitting scheme, which solves either the (multi-phase) fluid flow or the mechanics problem first, and then updates the remaining solution variable(s) using the intermediate information \citep{prevost1997partitioned, settari1998coupled, minkoff2003coupled, suh2021asynchronous}. 
At a single time step, this procedure can be repeated iteratively until convergence is reached, or one may employ a staggered scheme without any iterations by treating the system of equations in a completely decoupled manner. 
The major upshot of this approach is that it can leverage existing simulation codes, only requiring the implementation of the interface between the two to ensure proper feedback and communication among the solvers. 
Adopting sequential schemes, however, is often limited by their inherent convergence issues and conditional stability, which can impact the overall reliability and robustness of the solution, particularly for strongly coupled systems \citep{li2003iterative, dean2006comparison, kim2011stability}. 
Another strategy is to use a monolithic scheme that requires a unified solver that updates solution variables simultaneously at each time step. 
In comparison to the sequential solution strategy, the monolithic approach offers unconditional stability and convergence in time, provided the Ladyzhenskaya-Babu\v{s}ka-Brezzi (LBB) condition is satisfied \citep{babuvska1973finite, brezzi1985two, bathe2001inf}. 
Common choices to meet stability requirements include adopting stabilization techniques that introduce additional terms to eliminate spurious modes \citep{pastor1999stabilized, truty2006stabilized, white2008stabilized}, or using inf-sup stable finite element discretizations, such as Taylor-Hood or Raviart-Thomas elements \citep{murad1994stability, borja1998elastoplastic, barbeiro2010priori, suh2021immersed}. 
However, the monolithic approach may involve substantial programming efforts, as stabilization techniques often lead to a modified variational formulation, while implementing LBB-stable elements requires specialized meshing and data structures to accommodate the different basis functions needed for the solution variables.

As an alternative, exploiting pre-trained large language models (LLMs)--deep learning models designed for natural language processing tasks--for numerical implementations has recently gained popularity as they demonstrated their capabilities in generating computer codes automatically from the provided set of natural language prompts \citep{bareiss2022code, buscemi2023comparative, chen2023improving, hou2023large, gu2024effectiveness}. 
For instance, \citet{kashefi2023chatgpt} investigated the effectiveness of a pre-trained LLM in implementing numerical algorithms across different programming languages, and their study pointed out that the LLM may struggle or experience irregular interruptions, particularly when tasked with generating lengthy computer code scripts. 
\citet{kim2024chatgpt} demonstrated the feasibility of generating MATLAB codes for seepage analysis, slope stability assessments, and image segmentation by prompting the pre-trained LLM, but their investigation into numerically solving a partial differential equation was limited to the finite difference method. 
\citet{orlando2023assessing} examined the finite element implementation capabilities of the LLM for single-field problems by utilizing finite element libraries, such as deal.ii \citep{bangerth2007deal} and FEniCS \citep{logg2012automated}, and suggested that while the LLM can serve as an initial building block for finite element codes, it still requires human intervention and a solid understanding of the relevant domain knowledge. 
To address this issue, \citet{ni2024mechagents} presented a multi-agent modeling framework that organizes collaborative teams of artificial intelligence agents for elasticity problems, enabling them to write, execute, and self-correct finite element codes. 
Regardless, to the best of the authors' knowledge, no attempts have been made to explore the LLM's capability for numerically solving (hydro-mechanically) coupled multi-field problems that require the implementation of mixed finite elements with a monolithic solution scheme.

The goal of this paper is to explore the applicability of the pre-trained LLMs for implementing monolithic finite element models to simulate hydro-mechanically coupled processes, as an attempt to alleviate the difficulties associated with numerical implementations. 
Among a vast number of LLMs available, we choose ChatGPT o1 from OpenAI for finite element code generation, as it is known to be the most popular and is the latest model as of the date of this publication, despite the fact that some LLMs (e.g., StarCoder \citep{li2023starcoder} or GitHub Copilot \footnote{\url{https://copilot.github.com/}}) may outperform ChatGPT in this specific domain. 
In particular, this study explores how ChatGPT can be adopted to solve coupled initial boundary value problems (IBVPs) relevant to geotechnical engineering applications, particularly those related to stress-induced pore pressure dissipation and ground settlement. 
In the initial stage, information essential for the IBVP and finite element implementations is provided to ChatGPT for code generation. 
These information includes: (1) the strong form of the governing equation, (2) constitutive models, (3) domain geometry and material properties, (4) initial and boundary conditions, (5) spatiotemporal discretization, and (6) solution strategies. 
Then, we correct any identified errors or bugs by re-prompting ChatGPT with the initial error message generated from its output. 
This prompt augmentation process is repeated until the GPT-generated finite element model successfully passes the verification/validation test. 
By utilizing two distinct environments: FEniCS and MATLAB, we also explore the finite element code generating performance of ChatGPT with different levels of programming interfaces, as the former offers advanced features for finite element analysis while the latter does not.

The rest of the paper is organized as follows. 
In Section \ref{sec:formulation}, we provide a brief summary of the hydro-mechanically coupled formulation for unsaturated soils considered in this work, including the continuum representation, balance laws, and constitutive models. 
From this formulation, Section \ref{sec:model_problem_and_GPT_FEM} begins by presenting three model problems designed to test the code generation capability of ChatGPT, and then introduces the finite element implementation strategy through prompt augmentations. 
Then, Section \ref{sec:ver_val} presents the verification/validation test results for the finite element code generated by ChatGPT, followed by discussions on their performances and conclusions in Section \ref{sec:dis_and_con}.

As for notations and symbols, bold-faced letters denote tensors (including vectors which are rank-one tensors); the symbol `$\cdot$' denotes a single contraction of adjacent indices of two tensors (e.g.,\ $\vec{a} \cdot \vec{b} = a_{i}b_{i}$ or $\tensor{c} \cdot \tensor{d} = c_{ij}d_{jk}$); and the symbol `:' denotes a double contraction of adjacent indices of tensors of rank two or higher
(e.g.,\ $\tensor{C} : \vec{\varepsilon}$ = $C_{ijkl} \varepsilon_{kl}$). 
As for sign conventions, unless specified, the directions of the tensile stress and dilative pressure are considered positive.

\section{Hydro-mechanically coupled formulation for unsaturated soils}
\label{sec:formulation}
This section provides a brief overview of the hydro-mechanically coupled formulation for unsaturated soils considered in this work.  
We begin by outlining the continuum representation and necessary assumptions that form the basis for describing the motion of a three-phase porous medium within the framework of mixture theory (Section \ref{sec:cont_rep}). 
We then summarize the balance laws for its mass and linear momentum (Section \ref{sec:bal_laws}), followed by the material models for the phase constituents (Section \ref{sec:cons_eq}) based on an important premise in continuum-scale poromechanics: all mechanical processes are driven by the effective stress. 
The formulation presented in this section will lay the foundations for the model problems that will be presented later in Section \ref{sec:IBVP} of this study.

\subsection{Continuum representation}
\label{sec:cont_rep}
In this work, we consider unsaturated soil as a homogenized continuum mixture consisting of a solid matrix ($s$) with continuous pores filled with water ($w$) and air ($a$) phases. 
This allows us to define the volume fractions of each phase constituent as, 
\begin{equation}
\label{eq:vol_frac}
\phi^s = \frac{\mathrm{d} V_s}{\mathrm{d}V} 
\: \: ; \: \:
\phi^w = \frac{\mathrm{d} V_w}{\mathrm{d}V}
\: \: ; \: \:
\phi^a = \frac{\mathrm{d} V_a}{\mathrm{d}V},
\end{equation}
where $\mathrm{d}V = \mathrm{d}V_s + \mathrm{d}V_w + \mathrm{d}V_a$ is the total elementary volume of the mixture, which gives $\phi^s + \phi^w + \phi^a = 1$. 
We also define the saturation ratios of the in-pore phases (i.e., water and air) as,
\begin{equation}
\label{eq:sat}
S^w = \frac{\phi^w}{1 - \phi^s}
\: \: ; \: \:
S^a = \frac{\phi^a}{1 - \phi^s},
\end{equation}
such that $S^w + S^a = 1$. 
Further, by letting $\rho_s$, $\rho_w$, and $\rho_a$ denote the intrinsic mass densities of the solid, water, and air, respectively, the total mass density of the three-phase mixture $\rho$ can be defined as the sum of their partial densities (i.e., $\rho^\alpha = \phi^\alpha \rho_\alpha$ for $\alpha = \lbrace s, w, a \rbrace$):
\begin{equation}
\label{eq:mass_density}
\rho = \rho^s + \rho^w + \rho^a 
= \phi^s \rho_s + \phi^w \rho_w + \phi^a \rho_a
= \phi^s \rho_s + (1-\phi^s) S^w \rho_w + (1-\phi^s) S^a \rho_a.
\end{equation}

To arrive at the coupled formulation for unsaturated soils considered in this study, we make several assumptions at this stage. 
First, we restrict the deformation of the three-phase mixture to the geometrically linear and elastic regime. 
Second, we assume that the evolution of the geomaterial of interest can be described by the deformation of its solid skeleton, so that the motions of water and air phases can be described relative to that of the solid. 
Third, the intrinsic density of air is negligible (i.e., $\rho_a = 0$) and its intrinsic pressure $p_a$ is assumed to be in equilibrium with atmospheric pressure, or $p_a = 0$, such that there is no need to track the air phase.

\subsection{Mass and momentum conservation laws}
\label{sec:bal_laws}
In the absence of mass exchanges among the constituents, mass balance equations for the solid and water phases can be written as,
\begin{align}
\label{eq:bal_mass_s}
&\dot{\rho}^s + \rho^s \diver{\vec{v}}_s = 0, \\
\label{eq:bal_mass_w}
&\dot{\rho}^w + \rho^w \diver{\vec{v}}_s + \diver{(\rho^w \tilde{\vec{v}}_w)} = 0, 
\end{align}
where $\dot{(\circ)} = \mathrm{D}^s ( \circ )/\mathrm{D} t$ is the material time derivative following the solid motion, and $\tilde{\vec{v}}_w = \vec{v}_w - \vec{v}_s$ is the relative velocity of the water phase constituent with respect to the solid velocity $\vec{v}_s$. 
Considering the barotropic flows, Eqs.~\eqref{eq:bal_mass_s} and \eqref{eq:bal_mass_w} specialize into the following form \citep{borja2006mechanical,song2014mathematical}:
\begin{align}
\label{eq:bal_mass_s2}
&\dot{\phi}^s + \frac{\phi^s}{K_s} \dot{p}_s + \phi^s \diver{\vec{v}}_s = 0, \\
\label{eq:bal_mass_w2}
&\dot{\phi}^w + \frac{\phi^w}{K_w} \dot{p}_w + \phi^w \diver{\vec{v}}_s + \diver{\vec{w}_w} = 0.
\end{align}
Here, $K_s$ and $K_w$ are the intrinsic bulk moduli of the solid and water phase constituents, respectively, $p_s$ and $p_w$ are their intrinsic pressures, while $\vec{w}_w$ is the Darcy's velocity that satisfies $\diver{\vec{w}}_w = \phi^w \diver{\tilde{\vec{v}}}_w + (\grad{\rho}^w \cdot \tilde{\vec{v}}_w)/\rho_w$. 
From Eq.~\eqref{eq:bal_mass_s2}, notice that the time rate of change of the water volume fraction can be expanded as,
\begin{equation}
\label{eq:phi_w_dot}
\dot{\phi}^w = \dot{\overline{(1-\phi^s)S^w}} = -\dot{\phi}^s S^w + (1 - \phi^s) \dot{S}^w
= (1 - \phi^s) \dot{S}^w + \left( \frac{\phi^s}{K_s} \dot{p}_s + \phi^s \diver{\vec{v}}_s\right) S^w.
\end{equation}
Substituting Eq.~\eqref{eq:phi_w_dot} into Eq.~\eqref{eq:bal_mass_w2} yields the following mass balance equation for the mixture:
\begin{equation}
\label{eq:bal_mass}
(1-\phi^s) \dot{S}^w + \frac{S^w \phi^s}{K_s} \dot{p}_s + \frac{\phi^w}{K_w}\dot{p}_w + S^w \diver{\vec{v}}_s + \diver{\vec{w}_w} = 0.
\end{equation}

Let $\tensor{\sigma}^s$ and $\tensor{\sigma}^w$ denote the partial Cauchy stress tensors for the solid and water phase constituents, respectively. 
Then, from the passive gas assumption that we made, the total Cauchy stress $\tensor{\sigma}$ can be obtained from the sum \citep{prevost1980mechanics, borja2004cam}: 
\begin{equation}
\label{eq:total_cauchy}
\tensor{\sigma} = \tensor{\sigma}^s + \tensor{\sigma}^w,
\end{equation}
while the intrinsic pressures can be readily obtained as,
\begin{equation}
\label{eq:intrinsic_p}
p_s = -\frac{1}{3 \phi^s}\tr{(\tensor{\sigma}^s)}
\: \: ; \: \:
p_w = -\frac{1}{3 \phi^w}\tr{(\tensor{\sigma}^w)},
\end{equation}
such that the total mean pressure can be expressed as:
\begin{equation}
\label{eq:mean_p}
p = -\frac{1}{3}\tr{(\tensor{\sigma})} = \phi^s p_s + \phi^w p_w.
\end{equation}
Based on this setting, balance of linear momentum in the absence of inertial forces can be written as,
\begin{align}
\label{eq:bal_mom_s}
&\diver{\tensor{\sigma}}^s + \rho^s \vec{g} + \vec{h}^s = \vec{0}, \\
\label{eq:bal_mom_w}
&\diver{\tensor{\sigma}}^w + \rho^w \vec{g} + \vec{h}^w = \vec{0},
\end{align}
where $\vec{g}$ is the gravitational acceleration, while $\vec{h}^s$ and $\vec{h}^w$ are the drag force vectors on solid and water phases, respectively, that satisfy $\vec{h}^s + \vec{h}^w = \vec{0}$. 
Thus, the momentum balance equation for the entire mixture can be obtained from the summation of Eqs.~\eqref{eq:bal_mom_s} and \eqref{eq:bal_mom_w}:
\begin{equation}
\label{eq:bal_mom}
\diver{\tensor{\sigma}} + \rho \vec{g} = \vec{0}.
\end{equation}

\subsection{Constitutive equations}
\label{sec:cons_eq}
For unsaturated soils, a thermodynamically consistent Bishop-type effective stress $\tensor{\sigma}'$ may be written as \citep{borja2006mechanical, nuth2008effective, borja2009effective}, 
\begin{equation}
\label{eq:eff_str}
\tensor{\sigma}' = \tensor{\sigma} + \left( 1 - \frac{K}{K_s} \right) \bar{p} \tensor{I},
\end{equation}
where $K$ is the bulk modulus of the skeletal structure formed by the solid phase, $\bar{p} = S^w p_w + S^a p_a = S^w p_w$ is the pore pressure, and $\tensor{I}$ is the second-order identity tensor. 
In this case, linear elastic mechanical responses of the solid skeleton can be described by the following stress-strain relation: 
\begin{equation}
\label{eq:lin_elas}
\tensor{\sigma}' = \lambda \tr{(\tensor{\varepsilon})} \tensor{I} + 2 \mu \tensor{\varepsilon}.
\end{equation}
Here, $\lambda$ and $\mu$ are the Lam\'{e} constants of the solid skeleton, while $\tensor{\varepsilon} = (\grad{\vec{u}}_s + \transpose{\grad}{\vec{u}}_s)/2$ is the infinitesimal strain tensor, the symmetric part of the gradient of the solid displacement $\vec{u}_s$, where $\vec{v}_s = \dot{\vec{u}}_s$. 
The use of geometric linearity may further lead to the following relation for volumetric deformations \citep{suh2021asynchronous, heider2020phase}:
\begin{equation}
\label{eq:sol_vol_frac}
\phi^s = \phi^s_0 ( 1 - \diver{\vec{u}}_s ), 
\end{equation}
where $\phi^s_0$ is the reference solid volume fraction.

Focusing on the laminar porous medium water flow with low Reynolds number, this study employs the generalized Darcy's law that can be expressed as,
\begin{equation}
\label{eq:darcy}
\vec{w}_w = -\frac{k_rk}{\mu_w}(\grad{p}_w - \rho_w \vec{g}),
\end{equation}
where $k_r$ is the relative permeability, $k$ is the intrinsic permeability, and $\mu_w$ is the dynamic viscosity of water. 
To incorporate the water retention behavior of the mixture, we adopt the van Genuchten model \citep{van1980closed}, which is one of the most common choices in modeling unsaturated geomaterials, to relate the pore water pressure, relative permeability, and pore water pressure: 
\begin{align}
\label{eq:vG_swcc}
&S^w = S^w_\text{res} + \left( S^w_\text{max} - S^w_\text{res} \right) \left[ 1 + \left( \frac{-p_w}{\alpha_\text{vG}} \right)^{n_\text{vG}} \right]^{-m_\text{vG}}, \\
\label{eq:vG_rel_perm}
&k_r = \left( \frac{S^w - S^w_\text{res}}{S^w_\text{max} - S^w_\text{res}} \right)^{\frac{1}{2}} \left\lbrace 1 - \left[ 1 - \left( \frac{S^w - S^w_\text{res}}{S^w_\text{max} - S^w_\text{res}} \right)^{\frac{1}{m_\text{vG}}} \right]^{m_\text{vG}} \right\rbrace^2, 
\end{align}
where $S^w_\text{res}$ and $S^w_\text{max}$ are the residual and maximum water saturation ratios, respectively, while $\alpha_\text{vG}$, $n_\text{vG}$, and $m_\text{vG}$ are fitting parameters related to the shape of the water retention curve.

\section{GPT-based finite element implementation} 
\label{sec:model_problem_and_GPT_FEM}
From the coupled formulation that we described previously, Section \ref{sec:IBVP} presents three model problems considered in this work (more details can be found in Appendix \ref{app:gov_eqs}). 
In particular, we focus on specific IBVPs including the dissipation of pore water pressure via water mass diffusion, time-dependent differential settlement of a strip footing, and gravity-driven seepage, which are chosen to cover a range of simple to moderately complex problems relevant to geotechnical engineering applications. 
Then, Section \ref{sec:GPT-based_FEM} introduces the GPT-based finite element modeling framework, where we particularly focus on the prompt engineering strategies that leverage the automated programming capability of the pre-trained large language model. 
The finite element codes generated from the strategized prompts to solve the IBVPs numerically will be subjected to either verification or validation tests (Section \ref{sec:ver_val}), and our observations will also be discussed therein.

\subsection{Model problems}
\label{sec:IBVP}
To test the code generation capability of ChatGPT, this study considers three model problems that are derived from the coupled formulation summarized in Section \ref{sec:formulation}. 
As summarized in Table \ref{tab:model_probs}, each model problem is designed to highlight different aspects of hydro-mechanically coupled soil behavior under varying conditions, ranging from simple to moderately complex scenarios, in terms of finite element implementation. 
Further, to maintain consistency with previous studies (e.g., \citep{suh2021asynchronous, white2016block, na2017computational, suh2024data}) and for notational convenience, the following symbols or letters will be reserved for the model problems: 
\begin{equation}
\vec{u} \equiv \vec{u}_s
\: \: ; \: \:
\vec{w} \equiv \vec{w}_w
\: \: ; \: \:
\phi \equiv 1 - \phi^s
\: \: ; \: \:
B \equiv 1 - \frac{K}{K_s}
\: \: ; \: \:
M \equiv \left( \frac{B-\phi}{K_s} + \frac{\phi}{K_w} \right)^{-1}
\: \: ; \; \:
c_v \equiv K\frac{k}{\mu_w},
\end{equation}
such that $\phi$ indicates the porosity, $B$ the Biot's coefficient, $M$ the Biot's modulus, and $c_v$ the coefficient of consolidation \citep{biot1941general, terzaghi1943theoretical}. 

\begin{table*}[!t]%
\newcommand{\multrow}[1]{\begin{tabular}{@{}l@{}} #1 \end{tabular}}
\centering %
\caption{Three model problems considered in this study.\label{tab:model_probs}}%
\begin{tabular*}{\textwidth}{@{\extracolsep\fill}llll@{}}
\toprule
\scriptsize \textbf{Model problem} 
&\scriptsize \textbf{\multrow{(1) Fluid mass diffusion\\in 1D space}}  
&\scriptsize  \textbf{\multrow{(2) Hydro-mechanical coupling\\in saturated medium}}  
&\scriptsize  \textbf{\multrow{(3) Hydro-mechanical coupling\\in unsaturated medium}} 
\\
\midrule
\scriptsize \multrow{Formulation and\\phase constituents}
& \scriptsize Single-phase ($w$)
& \scriptsize Two-phase ($s$, $w$)
& \scriptsize Pseudo-three-phase ($s$, $w$, $a$) \\
\midrule
\scriptsize Prime variables 
& \scriptsize Pore water pressure ($p_w$)
& \scriptsize \multrow{Displacement ($\vec{u}$) and\\pore water pressure ($p_w$)}
& \scriptsize \multrow{Displacement ($\vec{u}$) and\\pore water pressure ($p_w$)} \\
\midrule
\scriptsize Governing equation(s)
& \scriptsize $
\dot{p}_w = c_v\frac{\partial^2 p_w}{\partial z^2}
$  
& \scriptsize $
\begin{cases}
\diver{(\tensor{\sigma}' - B p_w \tensor{I})} + \rho \vec{g} = \vec{0} \\
\frac{1}{M}\dot{p}_w + B \diver{\dot{\vec{u}}} + \diver{\vec{w}} = 0
\end{cases}
$
& \scriptsize $
\begin{cases}
\diver{(\tensor{\sigma}' - S^w p_w \tensor{I})} + \rho \vec{g} = \vec{0} \\
\phi \dot{S}^w + S^w \diver{\dot{\vec{u}}} + \diver{\vec{w}} = 0
\end{cases}
$ \\
\midrule
\scriptsize Constitutive law(s)
& \scriptsize -
& \scriptsize $
\begin{cases}
\tensor{\sigma}' = \lambda \tr{(\tensor{\varepsilon})} \tensor{I} + 2 \mu \tensor{\varepsilon} \\
\vec{w} = -\frac{k}{\mu_w}(\grad{p}_w - \rho_w \vec{g})
\end{cases}
$
& \scriptsize $
\begin{cases}
\tensor{\sigma}' = \lambda \tr{(\tensor{\varepsilon})} \tensor{I} + 2 \mu \tensor{\varepsilon} \\
\vec{w} = -\frac{k_r k}{\mu_w}(\grad{p}_w - \rho_w \vec{g}) \\
S^w = S^w_\text{res} + \left( S^w_\text{max} - S^w_\text{res} \right) \left[ 1 + \left( \frac{-p_w}{\alpha_\text{vG}} \right)^{n_\text{vG}} \right]^{-m_\text{vG}} \\
k_r = \left( \frac{S^w - S^w_\text{res}}{S^w_\text{max} - S^w_\text{res}} \right)^{\frac{1}{2}} \left\lbrace 1 - \left[ 1 - \left( \frac{S^w - S^w_\text{res}}{S^w_\text{max} - S^w_\text{res}} \right)^{\frac{1}{m_\text{vG}}} \right]^{m_\text{vG}} \right\rbrace^2
\end{cases}
$ \\
\midrule
\scriptsize Example IBVP
& \scriptsize \multrow{Diffusion of excess\\pore water pressure\\$[$Figure \ref{fig:prob1_geom}$]$}
& \scriptsize \multrow{Time-dependent differential\\settlement of a strip footing\\$[$Figure \ref{fig:prob2_geom}$]$}
& \scriptsize \multrow{Gravity-driven seepage\\$[$Figure \ref{fig:prob3_geom}$]$} \\
\bottomrule
\end{tabular*}
\end{table*}

\begin{figure}[htbp]
\centering
\subfigure[]{\label{fig:prob1_geom}
\includegraphics[height=0.42\textwidth]{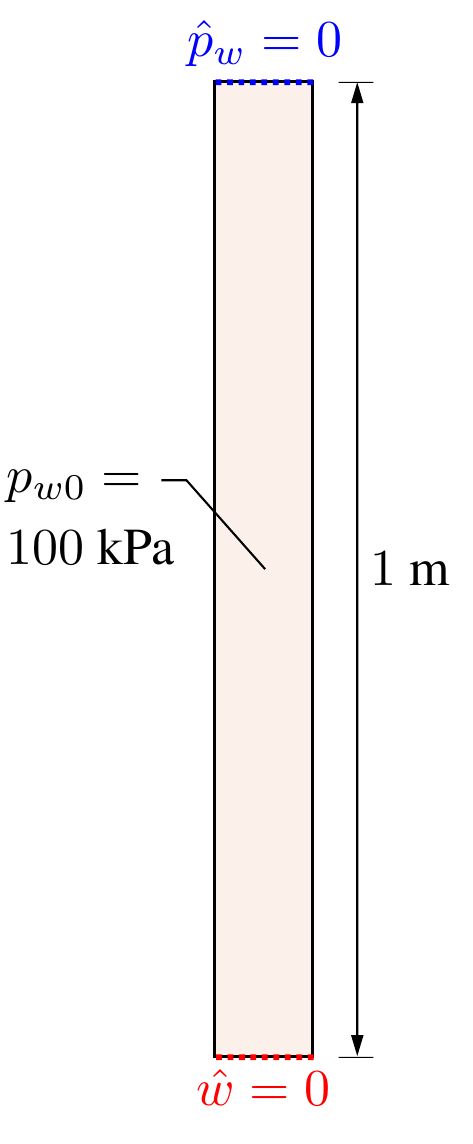}}
\hspace{0.01\textwidth}
\subfigure[]{\label{fig:prob2_geom}
\includegraphics[height=0.42\textwidth]{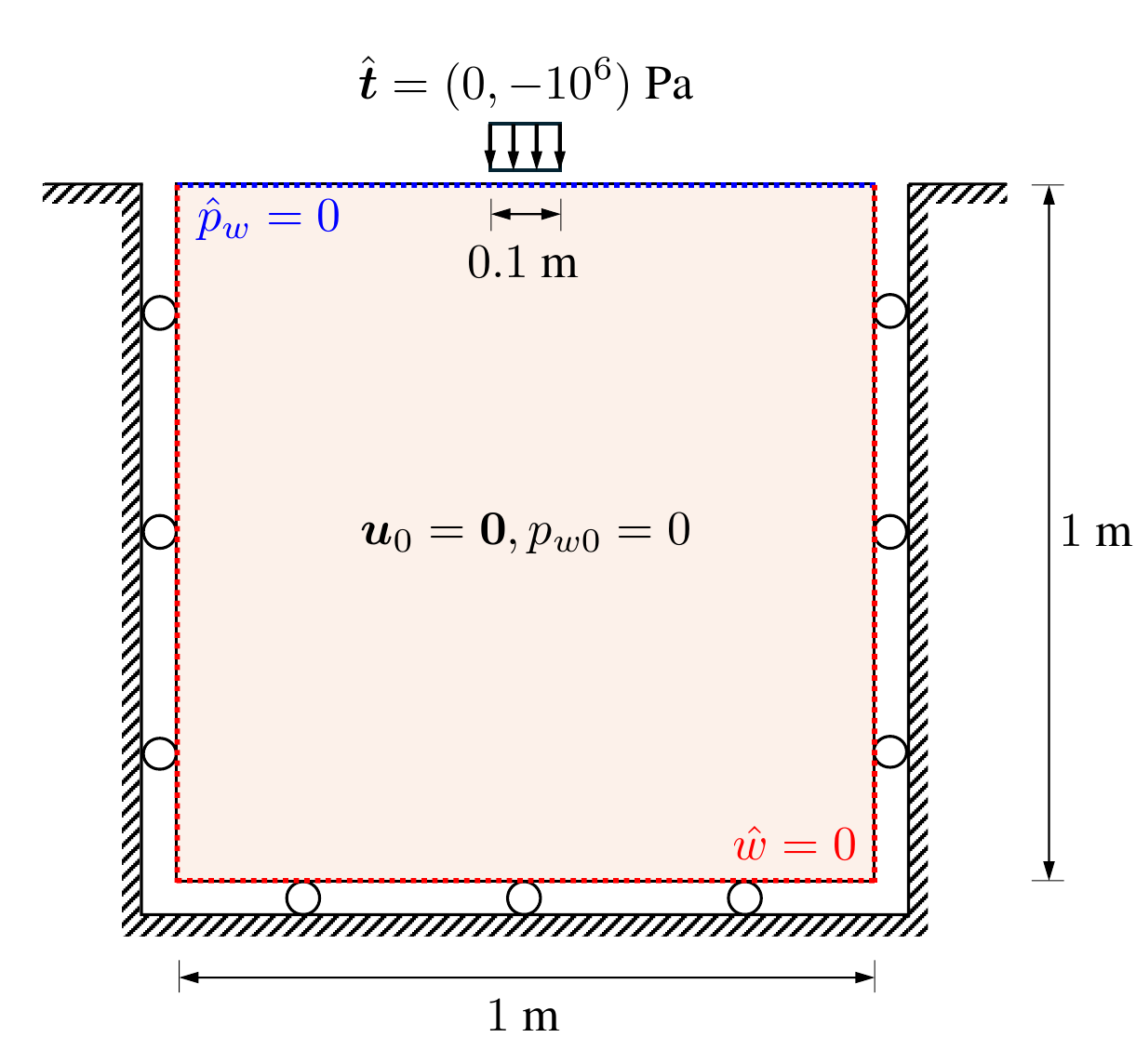}}
\hspace{0.01\textwidth}
\subfigure[]{\label{fig:prob3_geom}
\includegraphics[height=0.42\textwidth]{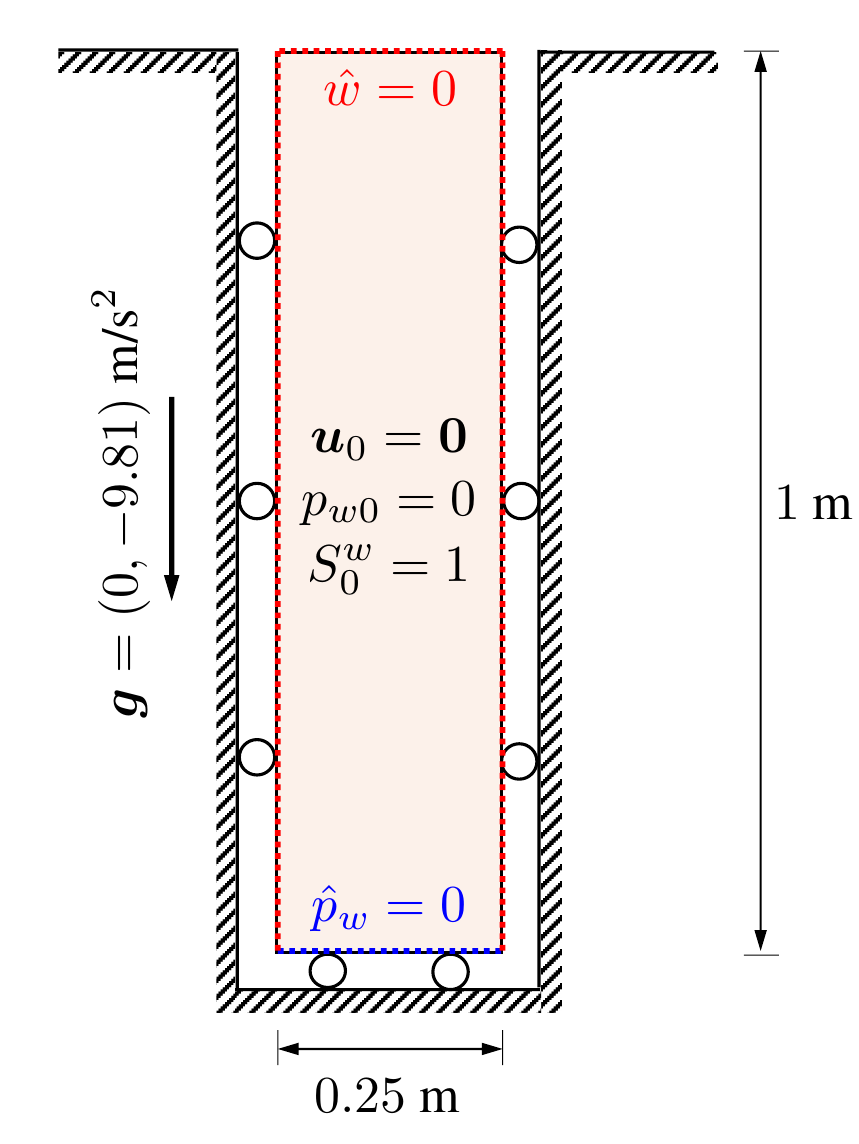}}
\caption{Schematic of geometry and boundary conditions for the numerical examples considered in this work: (a) dissipation of excess pore water pressure via fluid mass diffusion; (b) time-dependent differential settlement of a strip footing; and (c) gravity-driven seepage.}
\label{fig:geom}
\end{figure}

\textbf{Model Problem (1)} is the simplest among the problems considered in this work, which focuses on a fully saturated soil composed of incompressible phase constituents. 
In this case, the effective stress principle in Eq.~\eqref{eq:eff_str}, along with the linear elastic assumption in Eq.~\eqref{eq:lin_elas}, enables us to express the volumetric strain rate of the solid skeleton $\tr{(\dot{\vec{\varepsilon}})}$ in terms of the rate of change of pore water pressure $\dot{p}_w$. 
As outlined in Appendix \ref{app:model_prob1}, by limiting our attention to the volumetric deformation of the soil body, this integration allows us to simplify the coupled formulation into a single partial differential equation that considers $p_w$ as the prime variable, eliminating the need for a mixed finite element implementation. 
If we further limit the flow of pore water in the direction of loading only and neglect the gravitational effects, the governing partial differential equation reduces to a form identical to Terzaghi's 1D consolidation equation \citep{terzaghi1943theoretical}. 
From this simplified 1D mathematical model, we assign ChatGPT the task of writing a finite element code that can simulate the dissipation of excess pore water pressure over time in a meter-long soil column, as exemplified in Figure \ref{fig:prob1_geom}. 
Specifically, while prescribing an initial pore water pressure $p_{w0} = 100$ kPa, we replicate a single-drained condition by imposing zero pore water pressure ($\hat{p}_w = 0$) at the top and a no-slip condition ($\hat{w} = 0$) at the bottom.

\textbf{Model Problem (2)} considers a deformable body of saturated soil with compressible phase constituents. 
Without further assumptions, as outlined in detail in Appendix \ref{app:model_prob2}, this setting results in a mathematical model consisting of two coupled partial differential equations, requiring a mixed finite element implementation to solve this problem numerically. 
Based on this framework, we focus on simulating the rate-dependent mechanical responses of the soil induced by the redistribution of pore water therein. 
The corresponding numerical example is, therefore, chosen to demonstrate ChatGPT's finite element implementation capability in a fully coupled manner (e.g., using monolithic solution strategy). 
As illustrated in Figure \ref{fig:prob2_geom}, the problem domain is a 1 m wide and 1 m deep square-shaped body of saturated soil, with a 0.1 m wide strip footing at the top center subjected to a compressive load of $\hat{t}_y = - 10^6$ Pa. 
The bottom and two lateral surfaces are set to be impermeable ($\hat{w} = 0$), while the top surface is allowed to drain freely during the simulation by imposing a zero pore water pressure boundary condition ($\hat{p}_w = 0$). 
Since fully saturated conditions do not require taking water retention behavior [i.e., Eqs.~\eqref{eq:vG_swcc} and \eqref{eq:vG_rel_perm}] into account, this problem can be seen as of intermediate difficulty compared to the other two examined in this work.

\textbf{Model Problem (3)} addresses the coupled behavior of unsaturated soil, representing the most general and complex case among the other problems considered in this study. 
As delineated in Appendix \ref{app:model_prob3}, this problem assumes that the phase constituents are incompressible and relies on a passive gas assumption that leads to a pseudo-three-phase formulation. 
This results in a mathematical model comprising two governing partial differential equations that describe the mass and momentum balances of the three-phase mixture, while necessitating an additional set of material models [e.g., Eqs.~\eqref{eq:vG_swcc} and \eqref{eq:vG_rel_perm}] to capture the unsaturated flow. 
Based on this setting, we consider a corresponding IBVP that closely resembles the experiment conducted by \citet{liakopoulos1964transient}. 
The problem domain is a 0.25 m wide and 1 m deep rectangular body of soil, as shown in Figure \ref{fig:prob3_geom}. 
The numerical simulation starts by assuming that the body of interest is initially equilibrated with hydrostatic pore water pressure, maintaining a fully saturated condition at $t = 0$. 
We then allow the water phase to escape from the bottom surface by prescribing a zero pore water pressure boundary condition ($\hat{p}_w = 0$), while imposing a zero flux boundary conditions ($\hat{w} = 0$) at all other boundary surfaces. 

\subsection{GPT-based finite element modeling framework}
\label{sec:GPT-based_FEM}
\begin{figure}[htbp]
\centering
\includegraphics[width=0.85\textwidth]{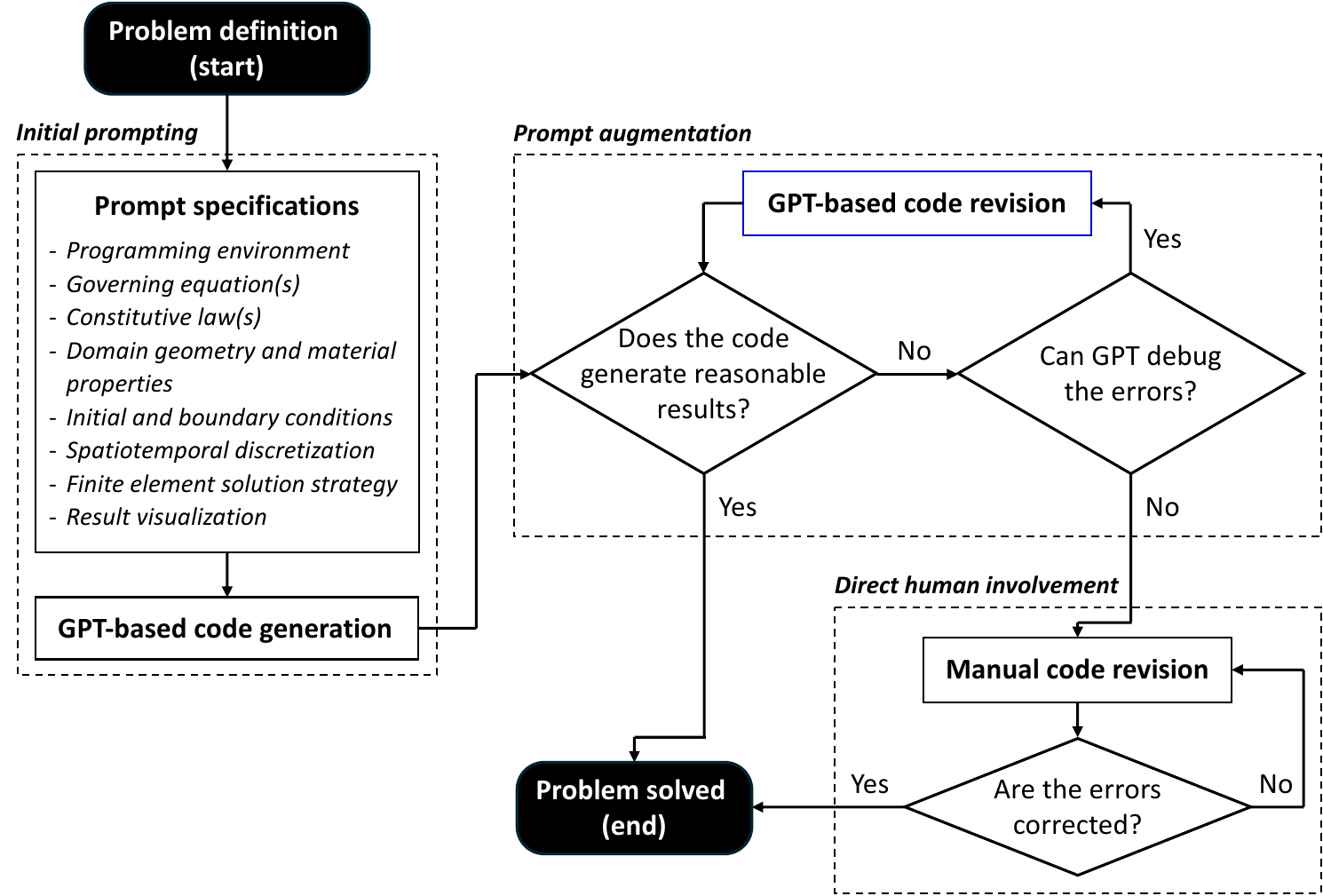}
\caption{Schematic representation of the workflow for GPT-based finite element implementation.}
\label{fig:gpt_flowchart}
\end{figure}

The ChatGPT-based finite element implementation strategy we use consists of three steps: initial prompting, prompt augmentation, and direct human involvement, indicated by dashed boxes in Figure \ref{fig:gpt_flowchart}. 
The initial prompting involves providing ChatGPT with essential information regarding the target IBVP and the corresponding solution strategy, while specifying the programming environment it should rely on. 
Then, we correct any errors or bugs in the initial code script generated by ChatGPT through a prompt augmentation process, using the error messages that appeared during execution. 
This ChatGPT-based debugging process is repeated until the finite element code passes the verification/validation tests. 
If the errors remain unresolved after several iterations, it indicates that direct intervention is required to fix the issues, which necessitates human programming efforts. 
The detailed procedure of the GPT-based finite element modeling framework (Figure \ref{fig:gpt_flowchart}) is described in the following.

\textbf{Initial prompting}. 
Firstly, we provide an initial prompt to task ChatGPT with implementing the finite element model using a specified programming environment. 
Key information that should be included in the prompt includes: (1) governing equation(s), (2) constitutive model(s), (3) domain geometry and material properties, (4) initial and boundary conditions, (5) spatiotemporal discretization, and (6) solution strategies. 
In particular, we provide the governing equations only in strong form, while the details on spatiotemporal discretization and solution strategies given to ChatGPT involve the type and size of the finite element, the time step size, the type of solver, and the time integration scheme. 
This specific setting has two purposes: one is to test whether ChatGPT can infer the weak form or the Galerkin form from the corresponding strong form, and the other is to introduce an implementation strategy targeted at users who are familiar with partial differential equations in their strong form but not necessarily with the variational formulations. 
After ChatGPT produces the initial code script from the given information, we execute the code to determine if any debugging is needed.

\textbf{Prompt augmentation}. 
If the finite element code generated from the initial prompt does not execute properly or generates inaccurate results, we then leverage ChatGPT's debugging capabilities via a series of prompt augmentation processes. 
This can be achieved by providing additional prompts that include either the error messages produced by the GPT-generated finite element code or the debugging instructions given by the user. 
This process is repeated iteratively until the code runs without any issues and yields results that are reasonably accurate. 
It should be noted that diagnosing errors may require an understanding of the formulation or the finite element method; however, this prompt augmentation process can eliminate the need for extensive programming efforts, while enabling users to explore the limits of ChatGPT's debugging capabilities.

\textbf{Direct human involvement}. 
This final step is necessary only if ChatGPT fails to resolve all issues after several iterations of the prompt augmentation (e.g., after 10 attempts). 
At this stage, we require direct human involvement in programming, following the same logic we apply to ChatGPT-based debugging. 
Nonetheless, this may demand less human effort than developing the finite element model manually from scratch, as the GPT-generated code can still serve as a useful starting point. 

\section{Verification and validation exercises}
\label{sec:ver_val}
This section presents the numerical results derived from GPT-generated codes for the model problems summarized in Section \ref{sec:IBVP} and discusses the hydro-mechanically coupled finite element modeling capabilities of ChatGPT.
Since the IBVPs considered herein have benchmark solutions or experimental results as counterparts, these problems serve as verification or validation exercises for the GPT-based finite element implementation framework. 
Specifically, we consider two different types of programming environments to test ChatGPT's adaptability: FEniCS and MATLAB, each having different levels of programming interfaces. 
As FEniCS offers advanced features such as the automatic evaluation of variational forms and the automated assembly of linear systems, it provides a high-level programming interface that simplifies scripting in Python. 
On the other hand, although MATLAB offers a wide range of pre-defined linear algebra tools and data structures, its functionality specific to finite element analysis is mostly low-level, requiring users to implement the building blocks of finite elements, from defining shape functions to assembling global matrices.

\subsection{Dissipation of excess pore water pressure via fluid mass diffusion}
\label{sec:prob1}
To generate a finite element code for solving Model Problem (1), we instructed ChatGPT to spatially discretize the geometry shown in Figure \ref{fig:prob1_geom} using 1D linear elements with a size of $\Delta z = 0.02$ m, and to use a time step size of $\Delta t = 1$ hour over a duration from $t_0 = 0$ to $t_f = 80$ hours, without providing any specific solution strategies.  
The only material property required for this simulation is the coefficient of consolidation, which was specified as $c_v = 0.016$ m$^2$/hr. 
Additionally, we asked ChatGPT to include a code block that generates the pore water pressure plots along the depth at $t = 5$, 10, 20, 40, and 80 hours. 

\begin{figure}[!htbp]
\centering
\includegraphics[width=1.0\textwidth]{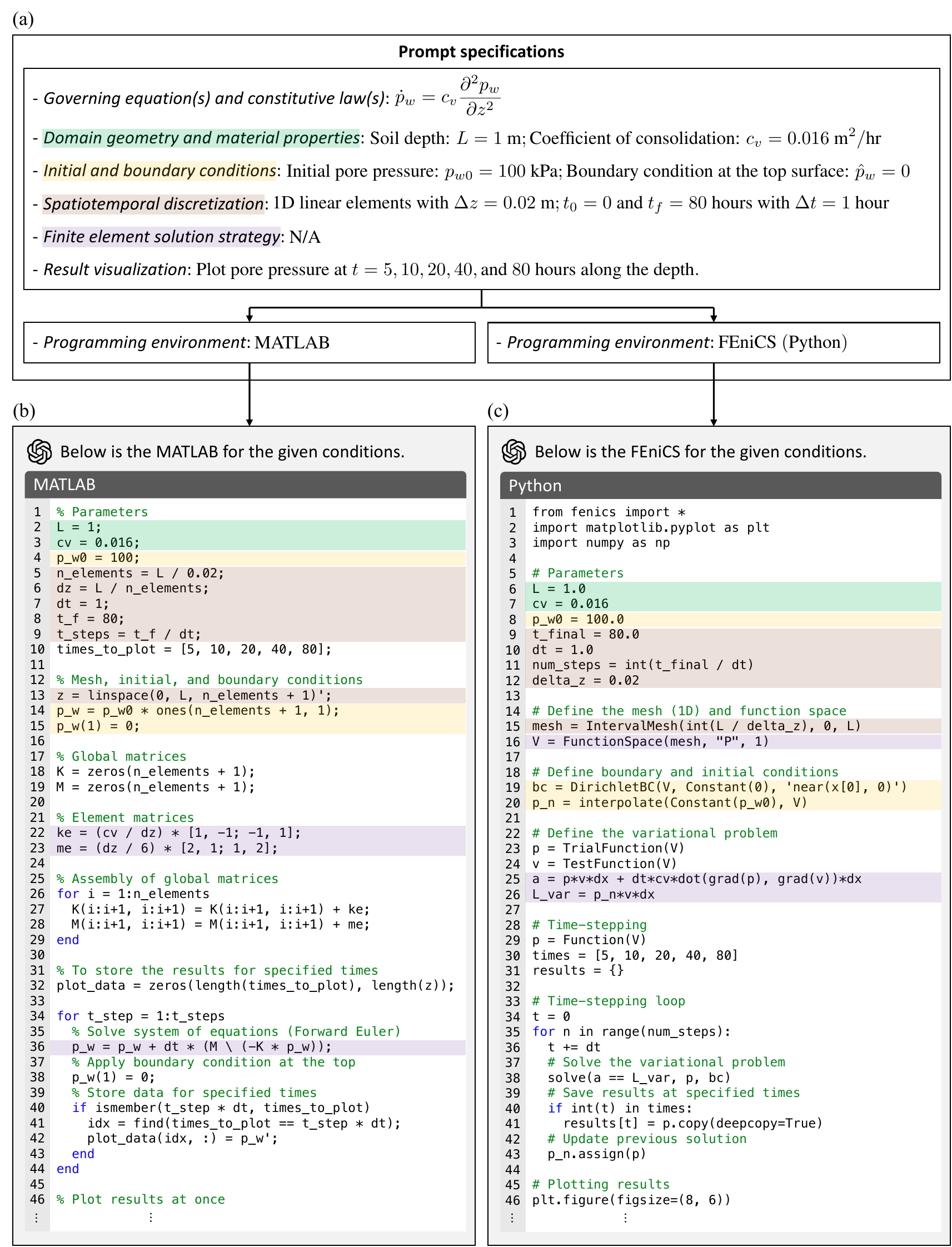}
\caption{GPT-based finite element code generation for Model Problem (1). (a) The initial prompts given to GPT; (b) an exemplary GPT-generated MATLAB script (GPT-MATLAB); and (c) the generated Python script using FEniCS (GPT-FEniCS).}
\label{fig:prob1_gpt_init}
\end{figure}

Figure \ref{fig:prob1_gpt_init} shows the specifications of the initial input prompt and the resulting MATLAB and FEniCS finite element code scripts that ChatGPT generated to numerically solve Model Problem (1). 
Here, the specific GPT-generated code blocks that correspond to the information we provided are color-coded for clarity. 
The lines of code highlighted with green boxes indicate the sections that define the domain geometry and material properties, the yellow boxes denote the parts that specify the initial and boundary conditions, the brown boxes are those related to the spatiotemporal discretization, while the purple boxes highlight the finite element solution strategy that ChatGPT employed. 
Both the GPT-generated MATLAB (GPT-MATLAB) and FEniCS (GPT-FEniCS) codes correctly implemented the specified requirements in a comprehensive manner, along with the relevant comments included. 
Even though the governing partial differential equation was presented in its strong form, GPT-MATLAB successfully constructed element matrices using a linear Lagrange shape function and assembled them into global matrices [e.g., Lines 22--29 in Figure \ref{fig:prob1_gpt_init}(b)]. 
This implies that ChatGPT can infer the appropriate Galerkin form from the given strong form equation. 
Similarly, when asked to generate the code using FEniCS, ChatGPT utilized the same linear element and implemented the weak form correctly [e.g., Lines 16 and 25--26 in Figure \ref{fig:prob1_gpt_init}(c)]. 
A notable observation is that both GPT-MATLAB and GPT-FEniCS utilized the trapezoidal rule for time integration through finite differences, without any specific instructions. 

\begin{figure}[htbp]
\centering
\includegraphics[width=1.0\textwidth]{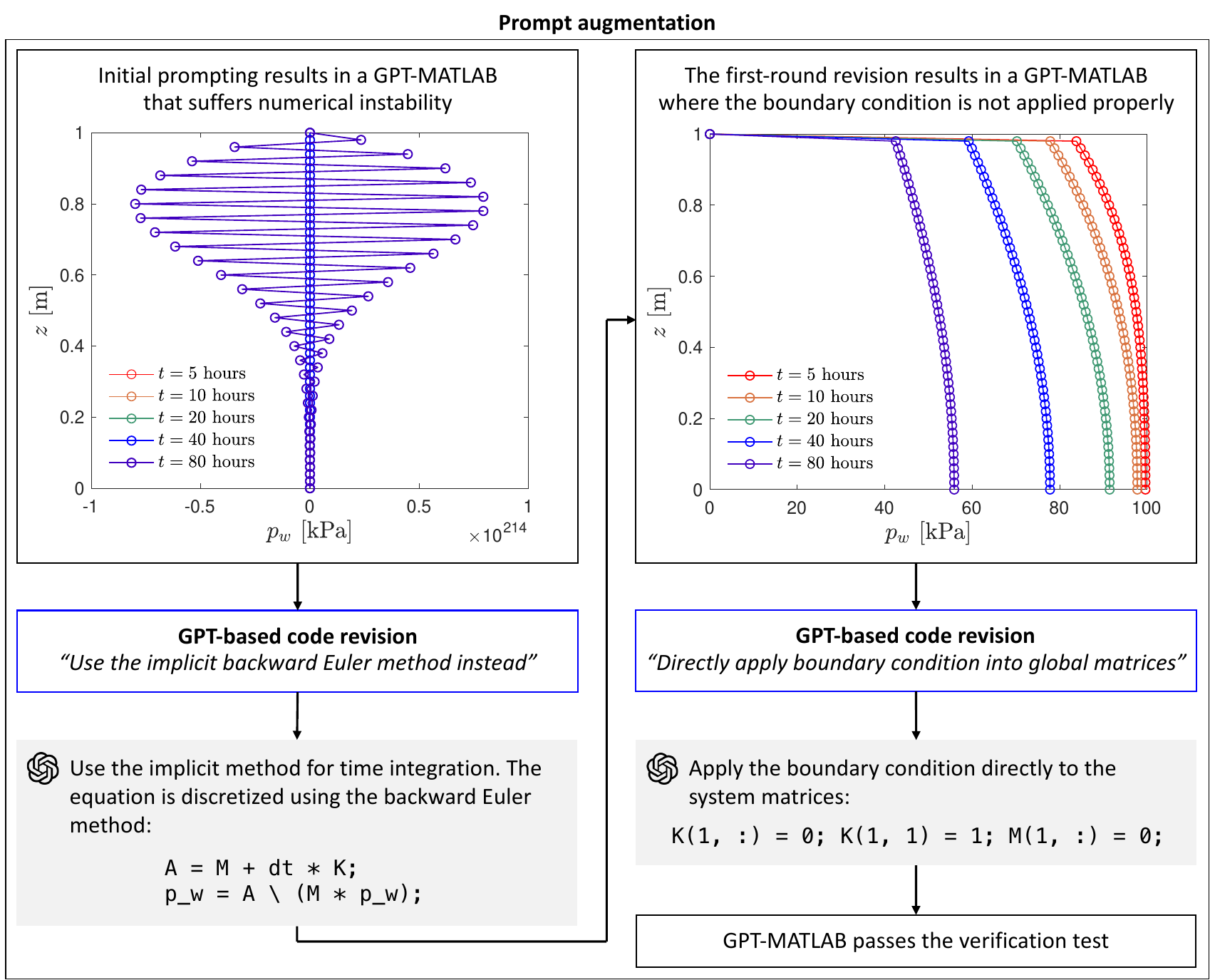}
\caption{Exemplary prompt augmentation process for GPT-MATLAB for Model Problem (1).}
\label{fig:prob1_gpt_aug}
\end{figure}

\begin{figure}[htbp]
\centering
\includegraphics[height=0.375\textwidth]{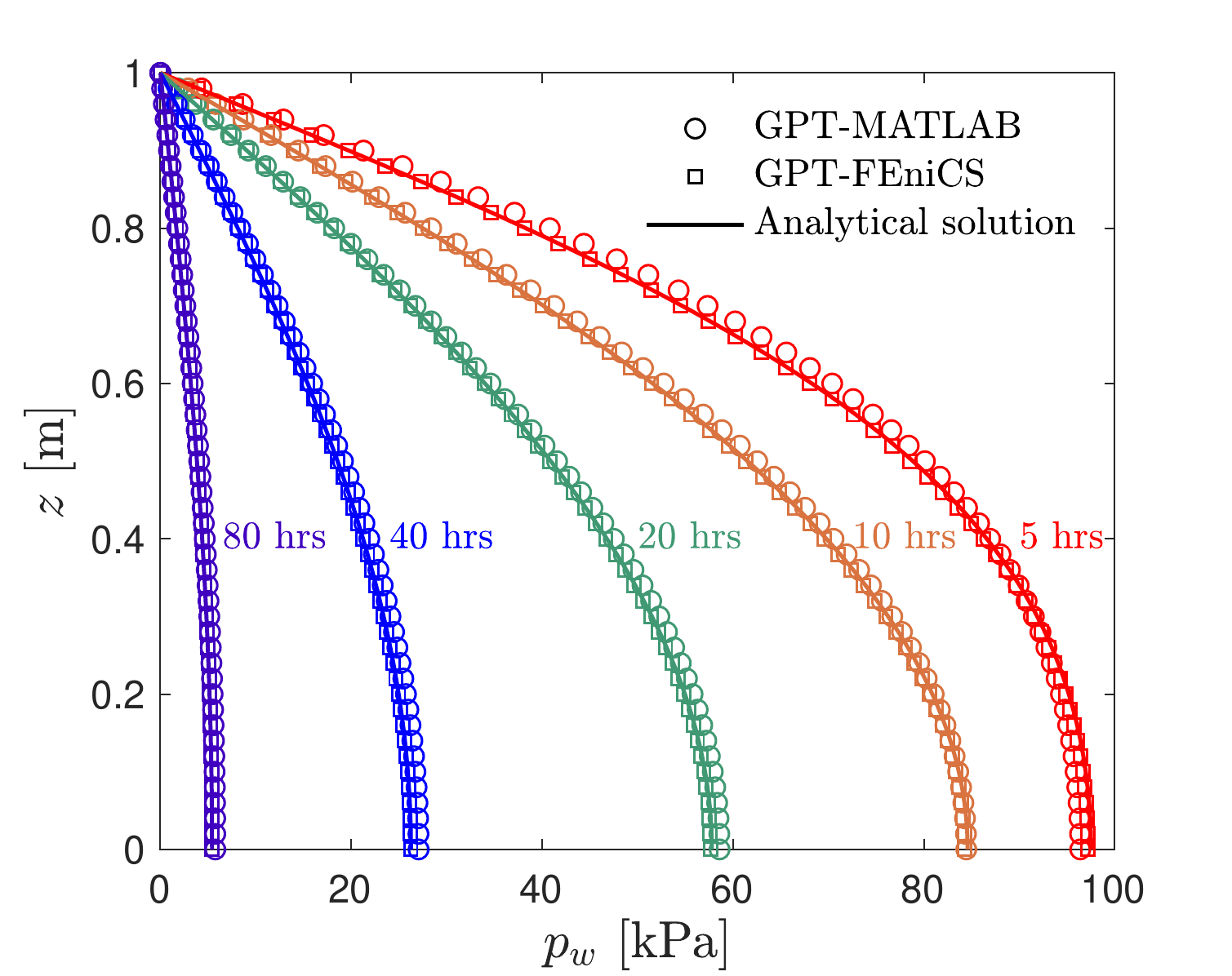}
\caption{Comparison of the final numerical results from GPT-MALAB and GPT-FEniCS against the analytical solution.}
\label{fig:prob1_final_result}
\end{figure}

While both GPT-written finite element codes in Figure \ref{fig:prob1_gpt_init} can be executed without any issues, GPT-MATLAB yielded unexpected results, whereas GPT-FEniCS passed the verification test on its initial attempt. 
To address the issues encountered by GPT-MATLAB, we proceeded to the second step, which involves a set of prompt augmentations as illustrated in Figure \ref{fig:prob1_gpt_aug}. 
Since the initial prompting results in a GPT-MATLAB code that experiences numerical instability due to the forward Euler scheme that it employed (i.e., the top-left box of Figure \ref{fig:prob1_gpt_aug}), we first provided an additional prompt to guide ChatGPT to use the implicit backward Euler method instead. 
As depicted in the top-right box of Figure \ref{fig:prob1_gpt_aug}, this first-round prompt augmentation step resolved the stability issues, though it did not address the issues caused by the boundary conditions that ChatGPT implemented improperly. 
We thus further instructed ChatGPT to apply the essential boundary condition directly to the pre-defined global matrices in this case, among multiple options. 
Consequently, ChatGPT replaced the equations corresponding to prescribed displacements in order to set the nodal pore water pressures to their prescribed values, thereby passing the verification test as illustrated in Figure \ref{fig:prob1_final_result}. 

\begin{figure}[htbp]
\centering
\includegraphics[height=0.375\textwidth]{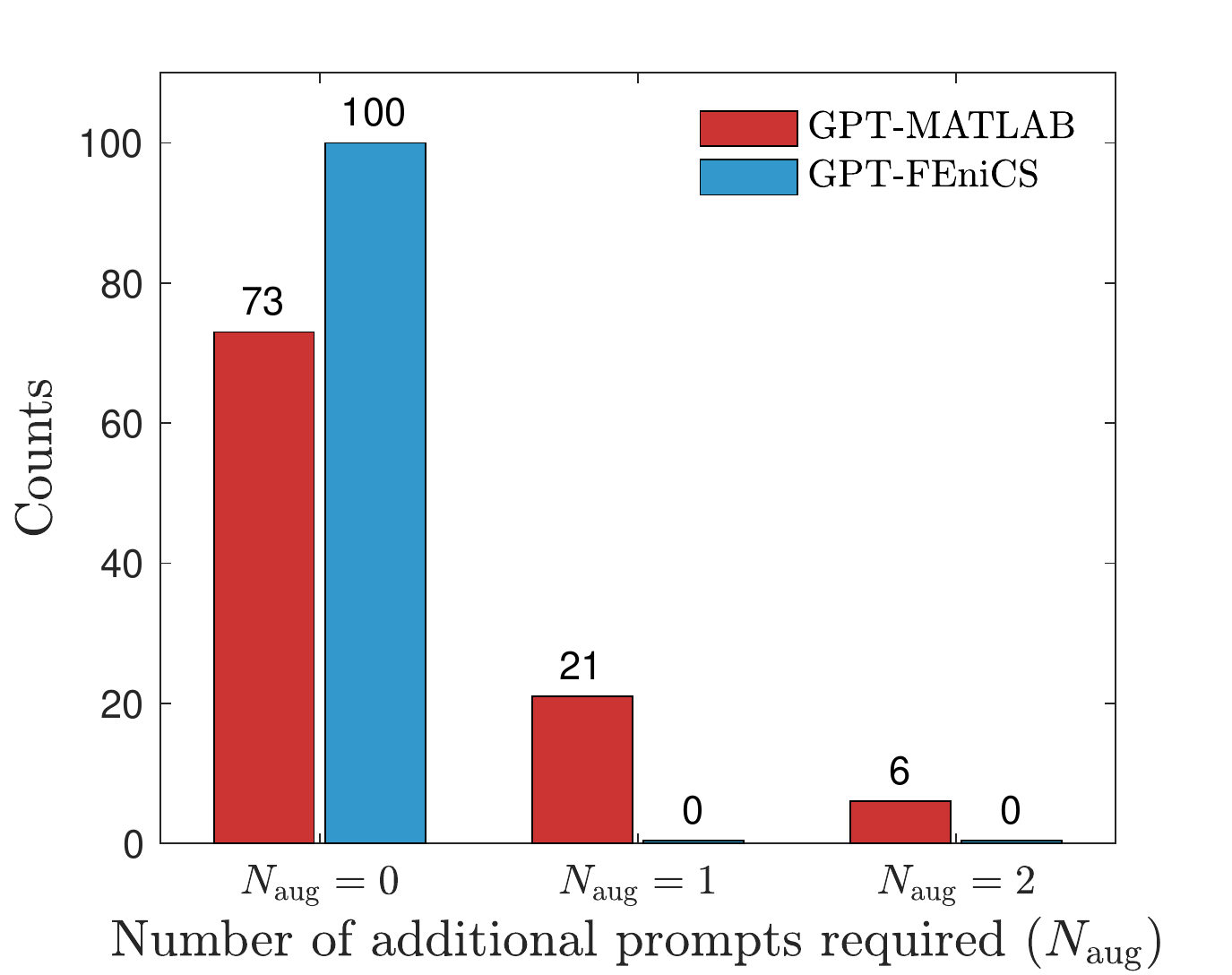}
\caption{Number of additional prompts required to pass the verification test for Model Problem (1).}
\label{fig:prob1_stats}
\end{figure}

It is important to note that ChatGPT may produce different finite element codes each time, as it is a probabilistic model. 
Hence, to test its consistency, we repeated the aforementioned processes 100 times in both the MATLAB and FEniCS environments. 
Figure \ref{fig:prob1_stats} shows the test result, summarizing the number of additional prompts ($N_\text{aug}$) required to pass the verification test. 
For GPT-MATLAB, 27 out of 100 cases required prompt augmentations to revise the code (i.e., $N_\text{aug} \ge 1$), whereas GPT-FEniCS passed the verification test on the first attempt (i.e., $N_\text{aug} = 0$) in all instances. 
While GPT-MATLAB may require further code revisions through re-prompting due to the low-level programming involved, the results suggest that ChatGPT can successfully implement a simple single-field finite element model, regardless of the programming environment used. 

\subsection{Time-dependent differential settlement of a strip footing}
\label{sec:prob2}
\begin{figure}[!htbp]
\centering
\includegraphics[width=1.0\textwidth]{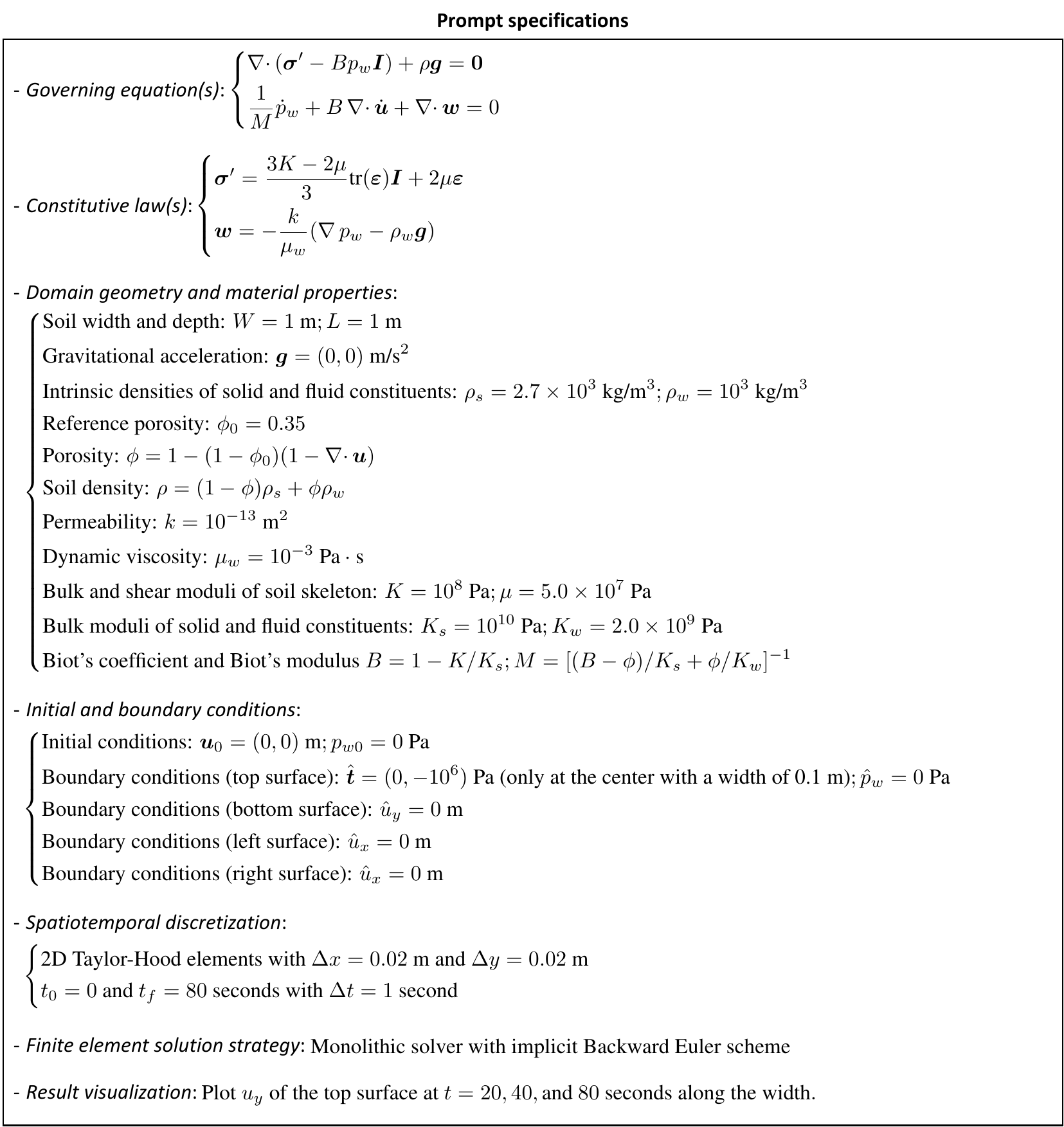}
\caption{The initial prompt given to GPT to generate mixed finite element code that can solve Model Problem (2).}
\label{fig:prob2_gpt_init}
\end{figure}

As Model Problem (2) involves two balance equations, two constitutive models, and a vast number of material properties, the initial prompt needed to guide ChatGPT in implementing the corresponding mixed finite element code is inevitably longer than that for the previous example, as illustrated in Figure \ref{fig:prob2_gpt_init}. 
For this model problem, we instructed ChatGPT to discretize the domain geometry [Figure \ref{fig:prob2_geom}] using a square-shaped, LBB-stable Taylor-Hood finite element with $\Delta x = \Delta y = 0.02$ m and to implement a time step size of $\Delta t = 1$ second over the duration from $t_0 = 0$ to $t_f = 80$ seconds. 
Furthermore, we explicitly directed ChatGPT to utilize an implicit backward Euler method to avoid the situation where it employs a time integration scheme that is only conditionally stable, while also asking it to generate a plot that shows the $y$-directional displacement of the top surface at $t = 20$, $40$, and $80$ seconds. 

\begin{figure}[htbp]
\centering
\includegraphics[height=0.375\textwidth]{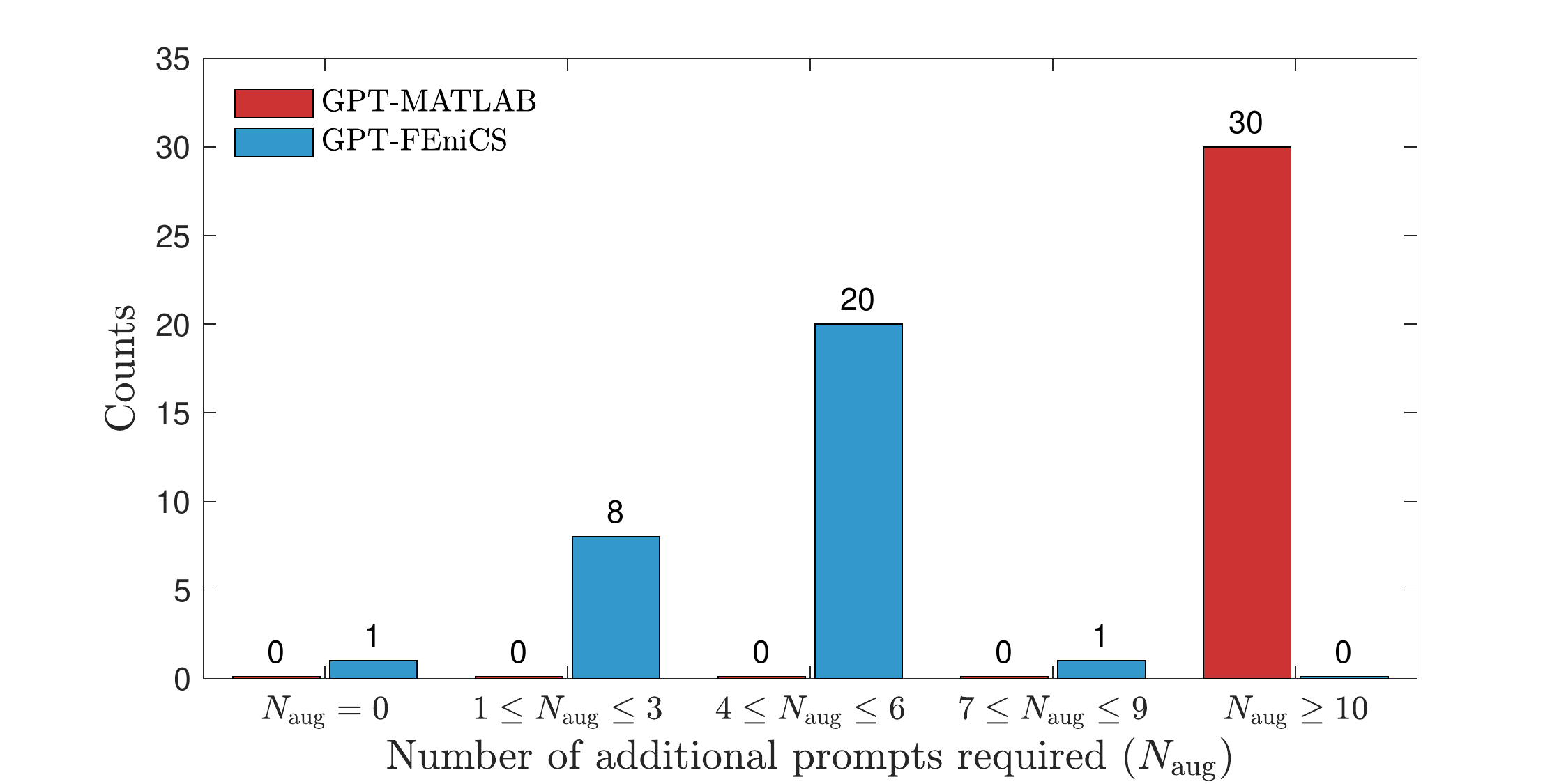}
\caption{Number of additional prompts required to pass the verification test for Model Problem (2).}
\label{fig:prob2_stats}
\end{figure}

Following the procedure outlined in Section \ref{sec:GPT-based_FEM}, the finite element implementation begins by providing ChatGPT with an initial prompt (Figure \ref{fig:prob2_gpt_init}), after which the generated code is revised through prompt augmentations for up to 10 iterations, until the code successfully passes the verification test. 
Similar to the previous exercise, this process was repeated 30 times in both the MATLAB and FEniCS environments to assess ChatGPT’s capability in generating mixed finite element code based on the programming interfaces used. 
Figure \ref{fig:prob2_stats} summarizes the number of additional prompts ($N_\text{aug}$) required in each case to successfully pass the verification test. 
In all 30 instances involving MATLAB, the results indicate that ChatGPT mostly fails to produce fully operational mixed finite element code, even after 10 rounds of prompt augmentations, demonstrating a consistent need for direct human intervention to address errors. 
Each GPT-MATLAB code script comprised over 300 lines, reflecting the inherent complexity associated with implementing mixed finite element models in this environment. 
This was primarily due to MATLAB's lack of features specifically designed for finite element analysis, which requires extensive low-level programming efforts, complicating the coding process and increasing the likelihood of errors. 
In contrast, the use of the FEniCS environment resulted in a marked enhancement of ChatGPT's performance. 
The FEniCS code generated by ChatGPT typically comprised fewer than 100 lines, demonstrating a prominent reduction in complexity compared to the GPT-MATLAB outputs. 
Moreover, GPT-FEniCS was able to pass verification tests with substantially fewer iterations, often achieving error-free codes within 6 prompt augmentations or fewer.
This improved efficiency can be attributed to FEniCS's higher-level abstractions and built-in features specialized for finite element analysis, which can simplify the coding process and consequently lower the risk of errors. 
The comparative analysis highlights the limitations of ChatGPT's code generation capabilities in low-level programming environments, emphasizing the advantages of leveraging existing domain-specific libraries. 

\begin{figure}[!htbp]
\centering
\includegraphics[width=1.03\textwidth]{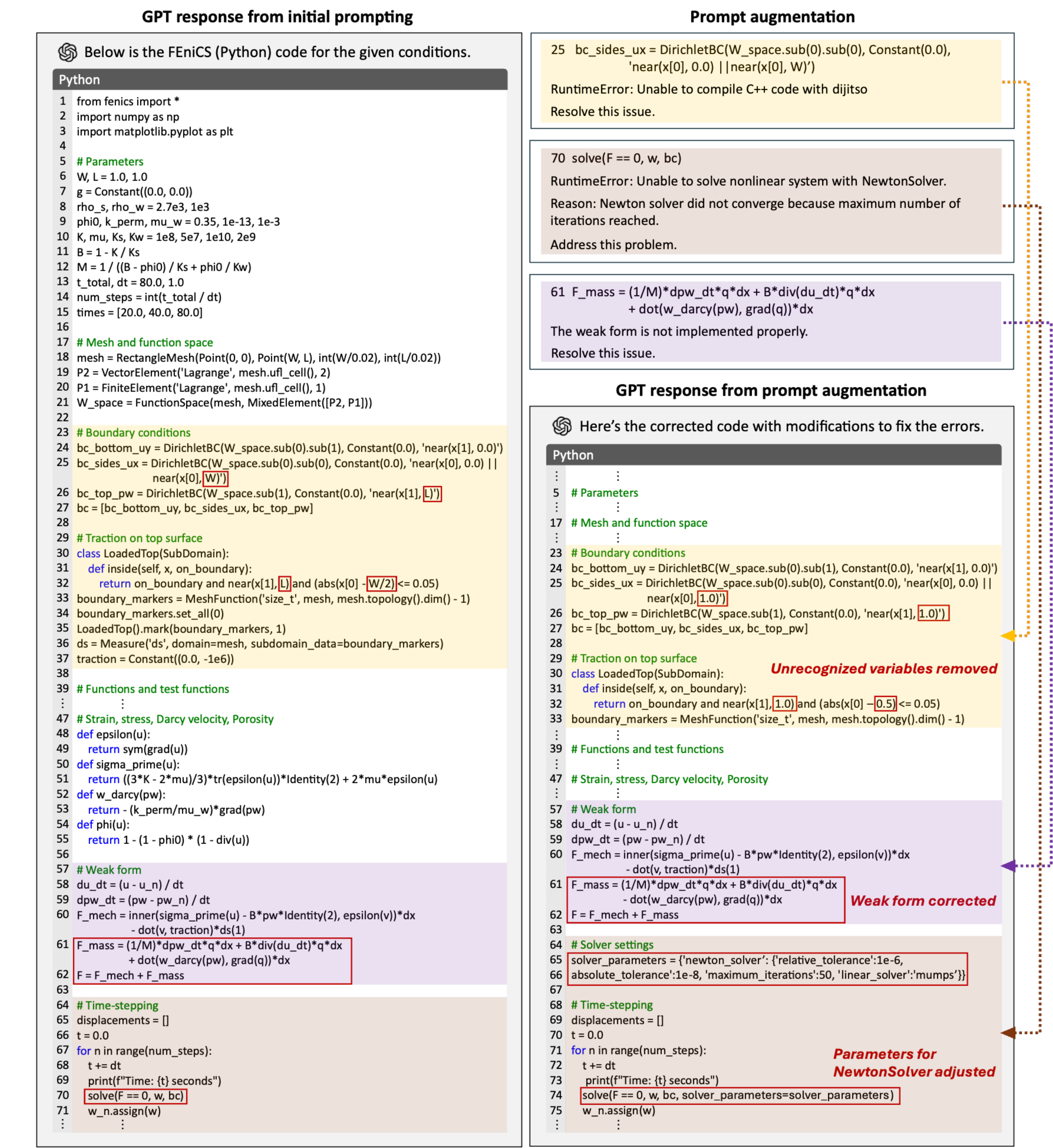}
\caption{Exemplary prompt augmentation process for GPT-FEniCS used to solve Model Problem (2).}
\label{fig:prob2_gpt_response}
\end{figure}

Figure \ref{fig:prob2_gpt_response} shows an exemplary GPT-based mixed finite element implementation procedure using FEniCS, which required three prompt augmentations to pass the verification test, with each augmentation clearly indicated by distinct color codes.  
At first, the GPT-generated code encountered execution failures due to a runtime error caused by using variables that cannot be recognized by the function \verb|DirichletBC| and the class \verb|LoadedTop| (e.g., red boxes within the yellow-colored region). 
To address this issue, we supplied ChatGPT with the specific line of code that caused the error, along with the generated error message, and requested to resolve the problem. 
A noteworthy observation is that, even though we provided ChatGPT with only line 25 of the code script depicted on the left side of Figure \ref{fig:prob2_gpt_response}, it replaced all the unrecognized variables (e.g., \verb|W| and \verb|L|), including those in lines 26 and 32. 
However, upon re-executing the revised GPT-FEniCS code, it encountered another runtime error because the \verb|NewtonSolver| utilized by ChatGPT failed to converge, attributed to the default tolerance being too low for this particular problem (e.g., red boxes within the brown-colored region). 
This required a second round of prompt augmentation, during which we again supplied ChatGPT with the relevant line of code and the error message it generated. 
In response to this feedback, ChatGPT specified the parameters for the solver, specifying lower values for \verb|relative_tolerance| and \verb|absolute_tolerance| compared to their default settings. 
These two rounds of re-prompting illustrate ChatGPT's capability to identify and resolve issues related to Python syntax, as well as its ability to utilize the features of the FEniCS library appropriately, highlighting its potential as a valuable tool in the development of mixed finite element code. 

\begin{figure}[htbp]
\centering
\subfigure[]{\label{fig:prob2_1st}
\includegraphics[height=0.375\textwidth]{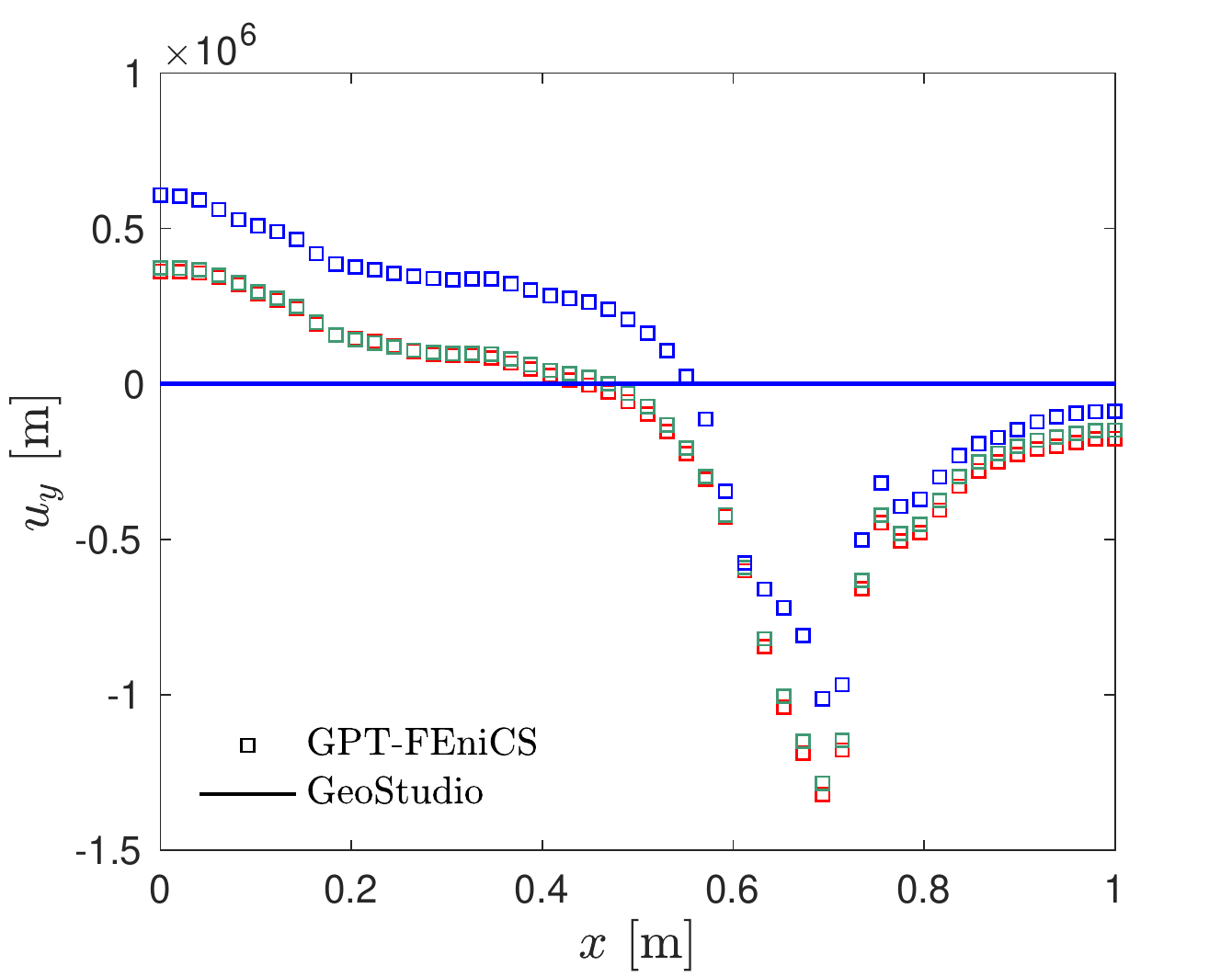}}
\hspace{0.01\textwidth}
\subfigure[]{\label{fig:prob2_2nd}
\includegraphics[height=0.375\textwidth]{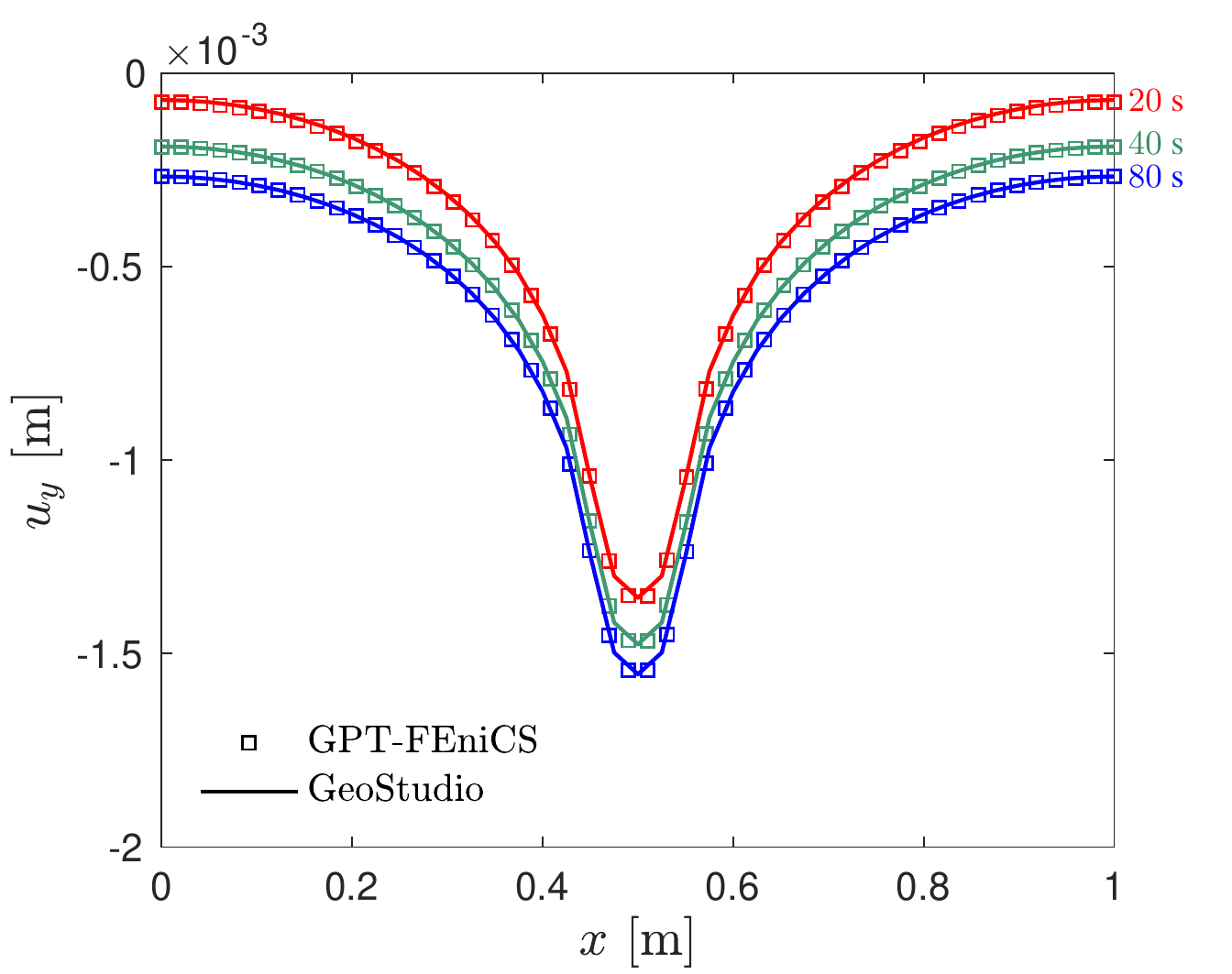}}
\caption{Numerical results from GPT-FEniCS after the (a) second and (b) third-round prompt augmentations, compared to the results from the commercial software GeoStudio.}
\label{fig:prob2_result}
\end{figure}

The two prompt augmentations resulted in a GPT-FEniCS code that executed without any issues; however, the results generated were significantly different from the benchmark solution obtained using the commercial software GeoStudio, as illustrated in Figure \ref{fig:prob2_1st}. 
This discrepancy was attributed to the incorrect implementation of the weak form by ChatGPT [cf. Eq.~\eqref{eq:prob2_weak}]. 
To resolve this issue, we instructed ChatGPT to correct this part by providing the specific line of code that required modification. 
As a result, the revised GPT-FEniCS code produced numerical results that closely aligned with the benchmark solution as illustrated in Figure \ref{fig:prob2_2nd}, successfully passing the verification exercise. 
This exercise suggests that while ChatGPT has significant potential in generating mixed finite element code by utilizing a high-level programming environment, it still requires careful oversight and iterative refinement to ensure both accuracy and reliability in its outputs. 

\subsection{Gravity-driven seepage}
\label{sec:prob3}
\begin{figure}[!htbp]
\centering
\includegraphics[width=1.0\textwidth]{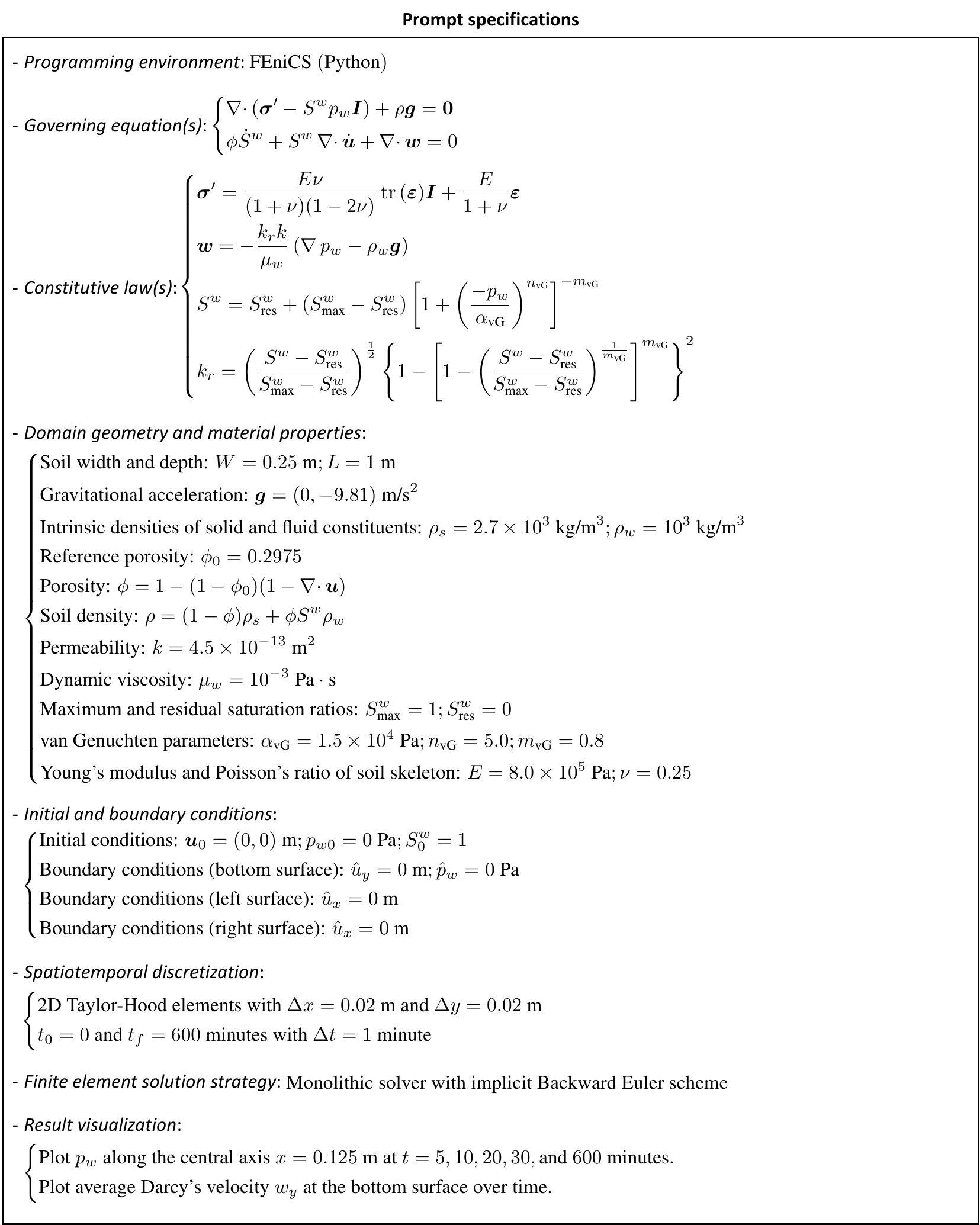}
\caption{The initial prompt given to GPT to generate mixed finite element code that can solve Model Problem (3).}
\label{fig:prob3_gpt_init}
\end{figure}

Based on the previous observations, this section focuses on generating a mixed finite element code for numerically solving Model Problem (3), solely relying on the FEniCS environment. 
As shown in Figure \ref{fig:prob3_gpt_init}, we directed ChatGPT to spatially discretize the geometry depicted in Figure \ref{fig:prob3_geom} using a Taylor-Hood finite element with $\Delta x = \Delta y = 0.02$ m and to employ a time step size of $\Delta t = 1$ minute for the simulation, covering the time interval from $t_0 = 0$ to $t_f = 600$ minutes. 
Similar to the previous exercise, we also instructed ChatGPT to solve the system monolithically while utilizing the backward Euler time integration scheme. 
Furthermore, we directed ChatGPT to generate a plot displaying the pore water pressure profile along the central axis (i.e., $p_w$--$y$ plot) at $t = 5$, $10$, $20$, $30$, and $600$ minutes, as well as a plot illustrating the average $y$-directional Darcy's velocity at the bottom surface over time (i.e., $w_y$--$t$ plot). 
Since this IBVP closely resembles the physical experiment conducted by \citet{liakopoulos1964transient}, the material parameters for the numerical simulation were chosen to align with the properties of the Del Monte sand as characterized their work, such that the simulation results can be compared with their measurements as a benchmark. 

\begin{figure}[htbp]
\centering
\includegraphics[height=0.375\textwidth]{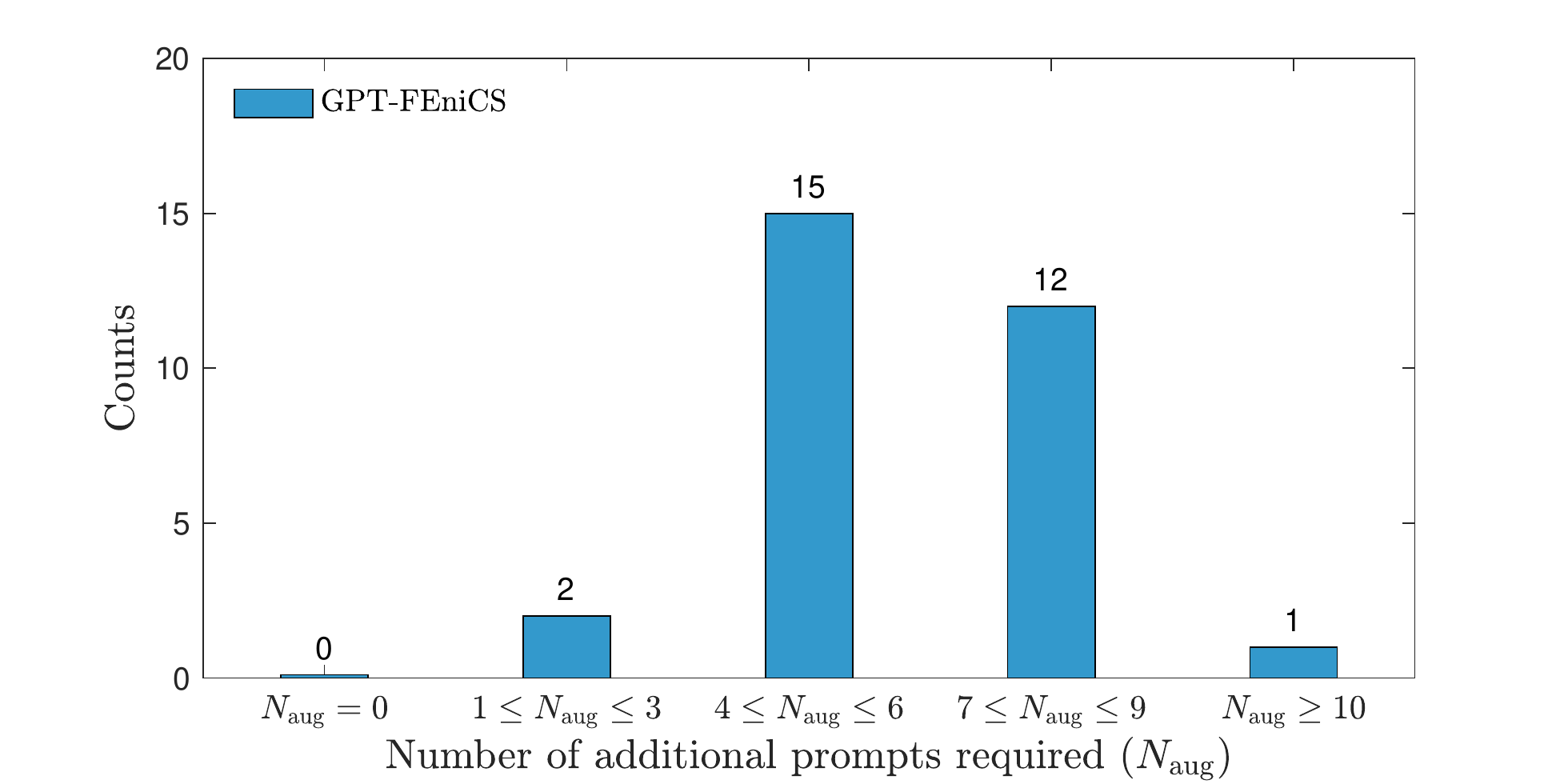}
\caption{Number of additional prompts required to pass the validation test for Model Problem (3).}
\label{fig:prob3_stats}
\end{figure}

Figure \ref{fig:prob3_stats} presents a bar plot illustrating $N_\text{aug}$ required to generate an error-free GPT-FEniCS code for numerically solving Model Problem (3), based on 30 repetitions of the initial prompting and subsequent augmentations. 
As this problem involves additional set of constitutive models to capture the water retention behavior of unsaturated soils, it required more re-prompting iterations compared to the previous exercise (e.g., Figure \ref{fig:prob2_stats}). 
The increased number of prompt augmentations suggests that even when leveraging a high-level programming interface like FEniCS, ChatGPT may struggle to generate error-free code on the first attempt, particularly when dealing with complex modeling requirements. 
Consequently, the trials exhibited a notable variation in $N_\text{aug}$ needed, ranging from as few as 1 to as many as 10, reflecting the complexities inherent in the problem being addressed. 
More importantly, none of the trials required direct human intervention, underscoring ChatGPT's capability to autonomously refine its outputs based on the error messages and specific guidance provided by the user during the re-prompting process. 
This result demonstrates a significant advancement in LLMs for numerical modeling, while also emphasizing the need for further improvements to enhance the reliability and efficiency of code generation for complex geotechnical engineering applications. 

\begin{figure}[!htbp]
\centering
\includegraphics[width=1.0\textwidth]{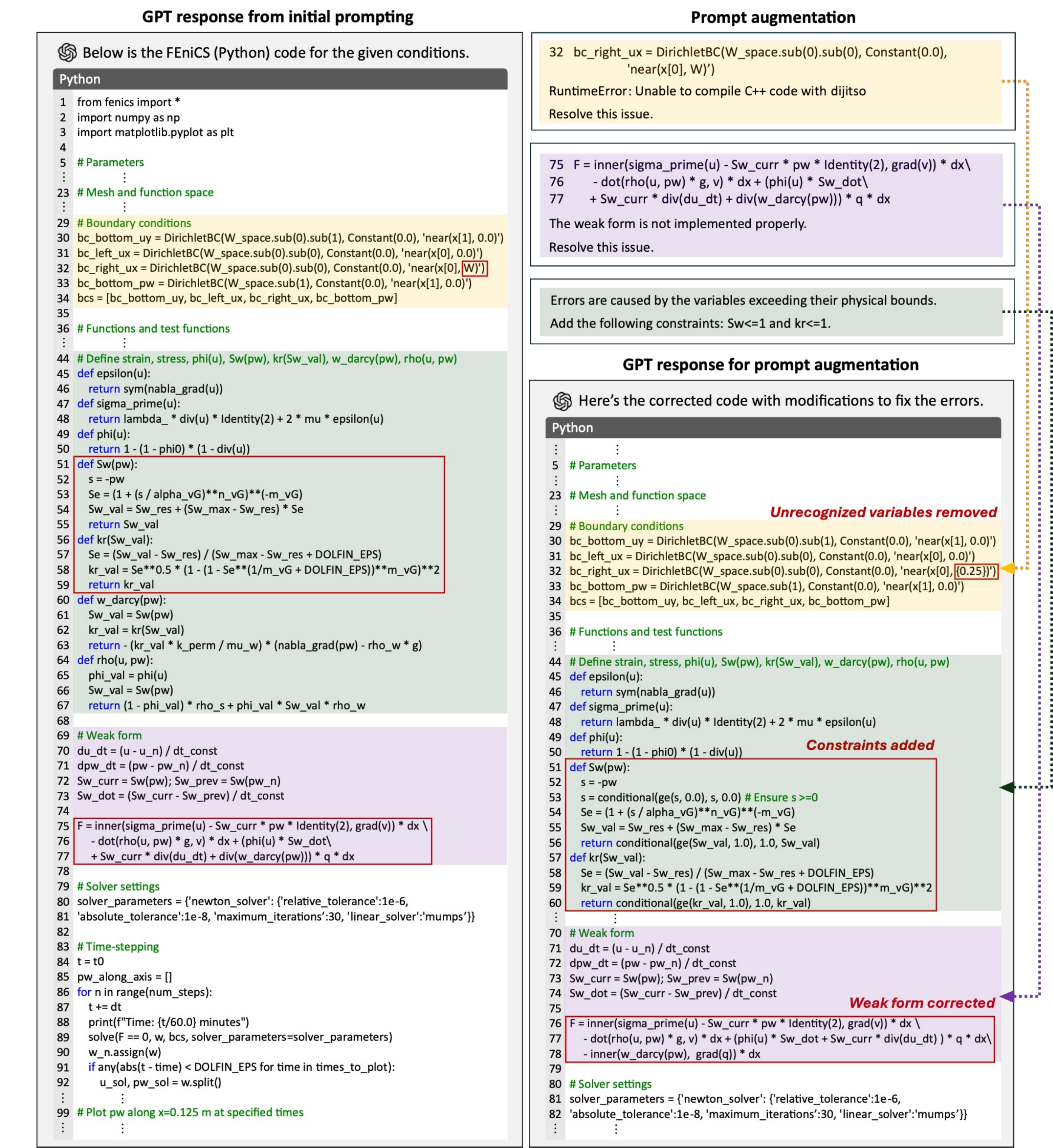}
\caption{Exemplary prompt augmentation process for GPT-FEniCS used to solve Model Problem (3).}
\label{fig:prob3_gpt_aug}
\end{figure}

\begin{figure}[htbp]
\centering
\subfigure[]{\label{fig:prob3_pw}
\includegraphics[height=0.375\textwidth]{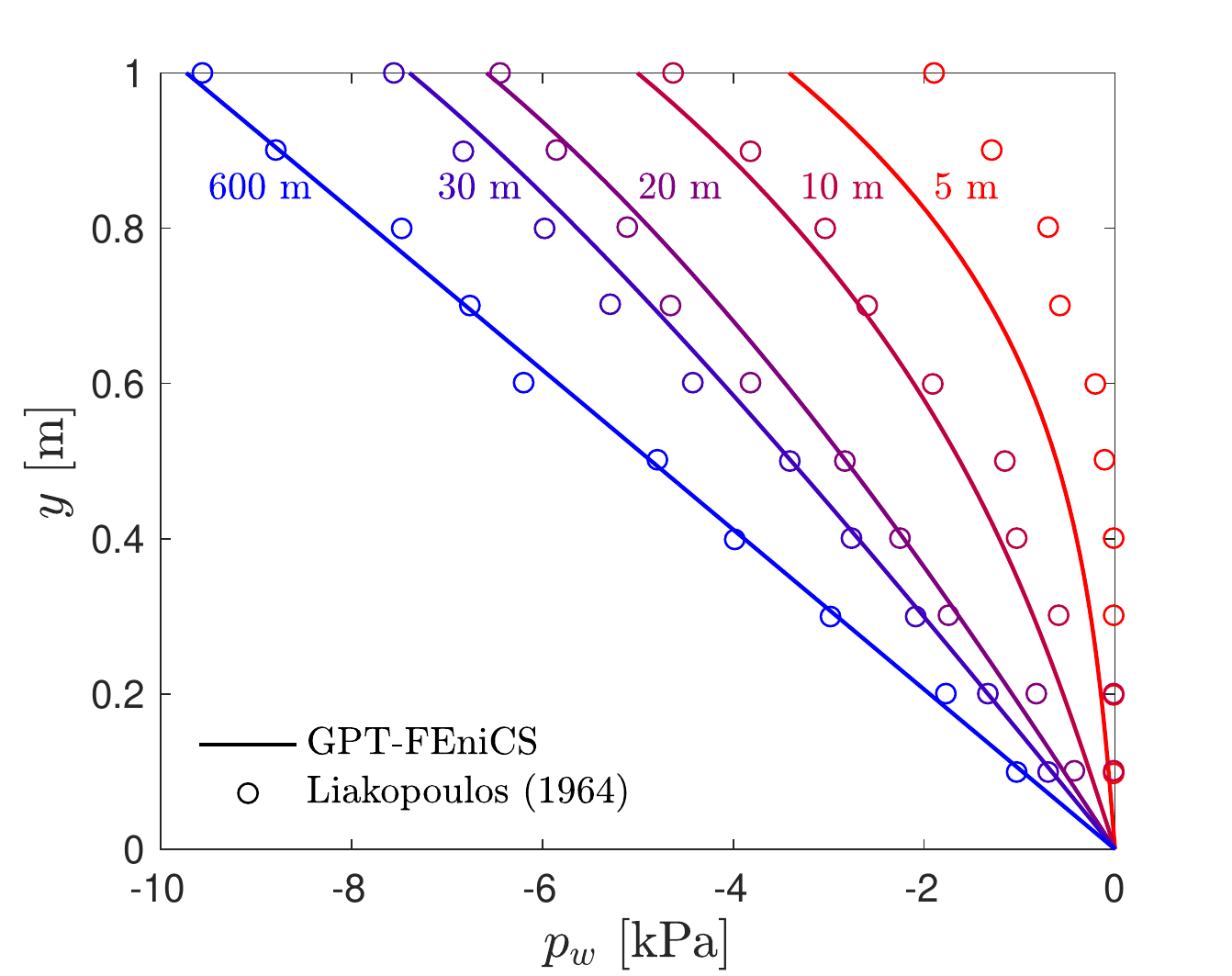}}
\hspace{0.01\textwidth}
\subfigure[]{\label{fig:prob3_wy}
\includegraphics[height=0.375\textwidth]{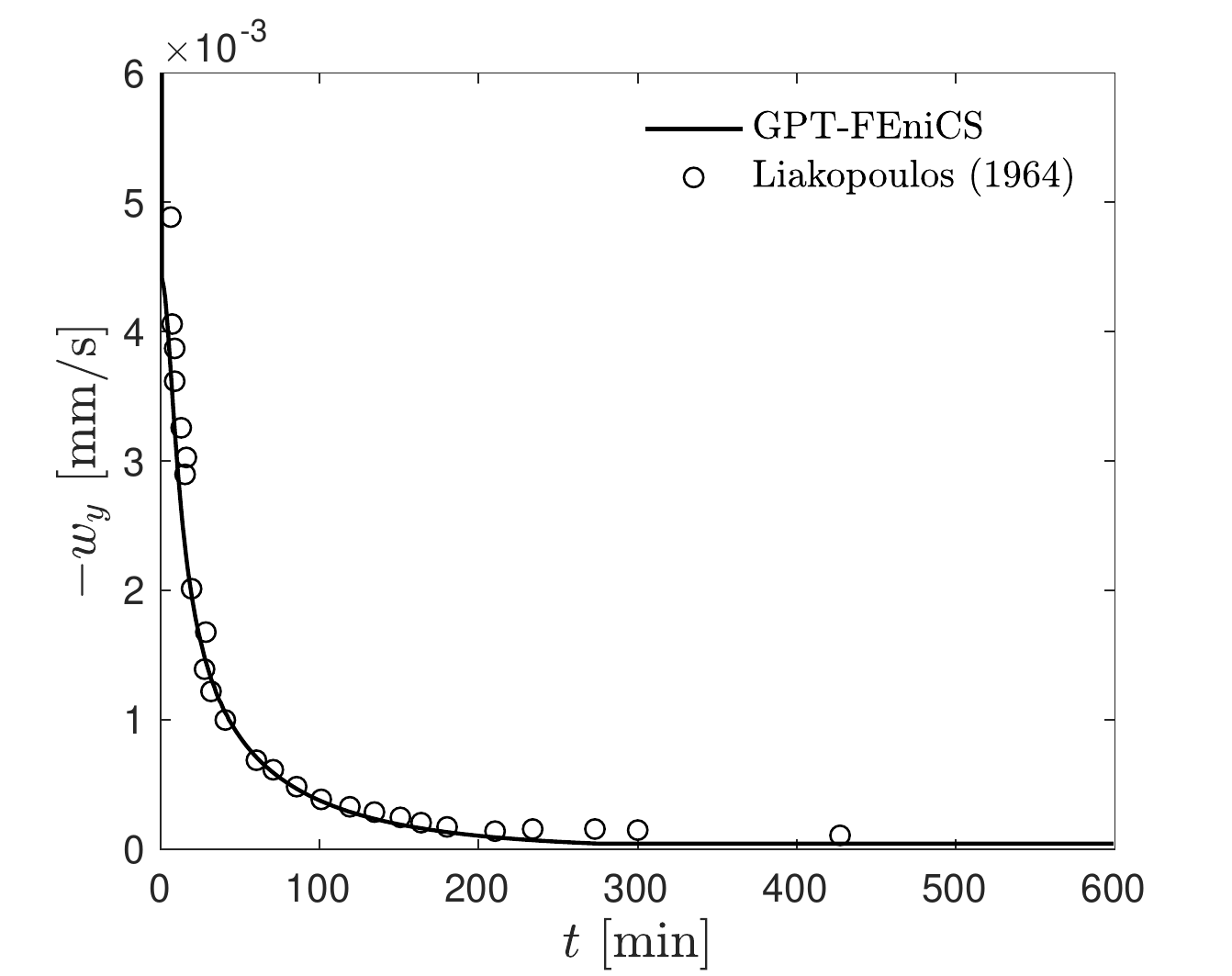}}
\caption{Numerical results from GPT-FEniCS after the three consecutive prompt augmentations, compared to the experimental results by \citet{liakopoulos1964transient}.}
\label{fig:prob3_result}
\end{figure}

Figure \ref{fig:prob3_gpt_aug} shows one of the trials conducted to generate a GPT-based mixed finite element code for Model Problem (3), which required three prompt augmentations to successfully pass the validation test. 
GPT-FEniCS code script generated after the initial prompting encountered a similar issue to that seen in the previous example (e.g., Figure \ref{fig:prob2_gpt_response}): a runtime error due to the use of a variable that cannot be recognized by the function \verb|DirichletBC|, as highlighted by the red box within the yellow-colored region in Figure \ref{fig:prob3_gpt_aug}. 
This issue was addressed in the same manner, by providing ChatGPT with the specific line of code responsible for the error (e.g., line 32) together with the corresponding error message it generated. 
However, the first-round prompt augmentation resulted in a GPT-FEniCS code that failed to execute. 
This was attributed to the omission of the integration by parts during the derivation of the weak form [cf. Eq.~\eqref{eq:prob3_weak}], as emphasized by the red box within the purple-colored region, which resulted in an insufficient representation of the specified natural boundary conditions and an inappropriate restriction of the solution space. 
In response to this issue, the second round of prompt augmentation focused on guiding ChatGPT to revise the weak form by specifying the lines of code that required modification (e.g., lines 75--77). 
Despite these two consecutive prompt augmentations, GPT-FEniCS was still not able to solve the system, due to the occurrence of the water saturation ratio exceeding the physical upper limit (i.e., $S^w > 1$), resulting in unphysical relative permeability values (i.e., $k_r > 1$). 
Hence, we instructed ChatGPT to impose the constraints $S^w \le 1$ and $k_r \le 1$ during the third round of prompt augmentation. 
This instruction led to the inclusion of conditional statements within the functions \verb|Sw(pw)| and \verb|kr(Sw_val)| that ChatGPT defined, which enforced that the degree of saturation and relative permeability derived from the van Genuchten model [i.e., Eqs.~\eqref{eq:vG_swcc} and \eqref{eq:vG_rel_perm}] remained within their permissible range of $[0, 1]$. 
As a result, the revised GPT-FEniCS code executed without any issues, producing numerical results that closely matched with the benchmark experimental data, as shown in Figure \ref{fig:prob3_result}, and successfully passing the validation test. 
This exemplary procedure suggests that, even though ChatGPT can autonomously generate a mixed finite element code with minimal or even no direct human programming efforts, it still necessitates appropriate guidance, which relies on the user's understanding of the underlying theoretical principles and the finite element method. 

\section{Discussion and conclusion}
\label{sec:dis_and_con}
This study investigated the capability of ChatGPT to generate (mixed) finite element code for solving hydro-mechanically coupled problems relevent to geotechnical engineering applications. 
We explored three model problems: one-dimensional consolidation, differential settlement of a footing, and gravity-driven unsaturated flow, utilizing both MATLAB and FEniCS programming environments. 
Our implementation strategy aimed to determine whether ChatGPT could produce reliable finite element code with minimal human intervention.

The findings indicate that ChatGPT can effectively generate complex mixed finite element codes when leveraging high-level programming interfaces, such as FEniCS. 
This capability is attributed to its inherent high-level abstraction and built-in functionalities, which significantly simplify the coding process. 
In contrast, when tasked with producing finite element code in the MATLAB environment, the GPT-generated code scripts often required extensive prompt augmentations and direct human involvement. 
This highlights the varying efficiencies of different programming environments and suggests that the choice of platform may significantly impact the success of automated finite element code generation.

Despite these promising results, the involvement of domain experts remains essential, particularly in diagnosing errors and instructing ChatGPT with necessary corrections. 
For complex problems, the limitations of LLMs become evident, as they may struggle to navigate intricate programming tasks without specified guidance. 
This dependence on human expertise emphasizes the need for collaborative approaches, where LLMs serve as assistive tools that complement the theoretical knowledge and practical experience of the users.

We believe that the proposed approach and the insights gained from this study can help remove technical or psychological barriers for readers interested in utilizing LLMs for finite element implementations. 
By highlighting the potential and limitations of the GPT-based finite element models, we aim to empower users to effectively integrate LLMs into their workflows.

In conclusion, while this study suggests that general-purpose LLMs like ChatGPT cannot yet fully replace human efforts in specialized domains, they can enhance efficiency in coupled finite element analysis. 
Effective utilization of these models may necessitate a synergistic approach, combining the capabilities of LLMs with the expertise of practitioners or researchers. 
As advancements in LLM technology continue, and as more sophisticated utilization strategies are developed, the gap between automated finite element code generation and expert-level programming is likely to narrow, paving the way for more integrated and efficient solutions in the field. 
Future research should focus on improving LLM training in specialized areas and exploring how these models can be better tailored to meet the specific needs of engineering applications.

\section*{Acknowledgments}
This work was supported by the National Research Foundation of Korea (NRF) grant funded by the Korea government (MSIT) (Nos. NRF-2021R1A5A1032433, 2023R1A2C2003534), and the start-up grant and Glennan Fellowship from Case Western Reserve University.

\bibliographystyle{unsrtnat}
\bibliography{main}

\begin{appendices}

\renewcommand{\theequation}{A\arabic{equation}}
\setcounter{equation}{0}
\renewcommand{\thefigure}{A\arabic{figure}}
\setcounter{figure}{0}
\section{Governing equations for model problems}
\label{app:gov_eqs}
\vspace*{12pt}
\subsection{Model Problem (1)}
\label{app:model_prob1}
The first model problem considered in this work assumes that the body of interest $\mathcal{B} \subset \mathbb{R}^3$ consists of incompressible phase constituents and is fully saturated with water, i.e., $K_s = K_w \approx \infty$ and $S^w = 1$. 
From these assumptions and by utilizing $\dot{\phi}^s + \dot{\phi}^w = 0$, the summation of Eqs.~\eqref{eq:bal_mass_s} and \eqref{eq:bal_mass_w} reduces the mass balance equations into the form identical to that seen in \citet{borja1995mathematical}:
\begin{equation}
\label{eq:prob1_eq1}
\diver{\dot{\vec{u}}} + \diver{\vec{w}} = 0.
\end{equation}
If we further neglect the gravitational effects (i.e., $\vec{g} = \vec{0}$), substituting Eq.~\eqref{eq:darcy} into Eq.~\eqref{eq:prob1_eq1} and utilizing $\diver{\dot{\vec{u}}} = \tr{(\dot{\tensor{\varepsilon}})}$ yields:
\begin{equation}
\label{eq:prob1_eq2}
\tr{(\dot{\vec{\varepsilon}})} = \frac{k}{\mu_w} \nabla^2{p}_w,
\end{equation}
where $\nabla^2 (\circ) = \nabla \cdot \nabla ( \circ )$ indicates the Laplacian operator. 
As Eq.~\eqref{eq:eff_str} reduces into $\tensor{\sigma}' = \tensor{\sigma} + p_w \tensor{I}$ in this case, notice that $\tr{(\dot{\tensor{\sigma}}')} = 3 \dot{p}_w$ under a constant total stress. 
This indicates that, from Eq.~\eqref{eq:lin_elas}, the volumetric strain rate of the solid skeleton and the pore water pressure can be related as,
\begin{equation}
\label{eq:prob1_eq3}
\tr{(\dot{\tensor{\varepsilon}})} = \frac{1}{3 \lambda + 2 \mu} \tr{(\dot{\tensor{\sigma}}')} = \frac{1}{K} \dot{p}_w,
\end{equation}
since $K = \lambda + 2\mu/3$. 
By substituting Eq.~\eqref{eq:prob1_eq3} into Eq.~\eqref{eq:prob1_eq2}, the two-phase problem reduces to a single-phase problem for the pore water pressure field $p_w$ as:
\begin{equation}
\label{eq:prob1_eq4}
\dot{p}_w = c_v \nabla^2 p_w.
\end{equation}
If we limit our attention to one-dimensional space (i.e., $\mathcal{B} \subset \mathbb{R}^1$), we arrive at the governing differential equation for the first model problem:
\begin{equation}
\label{eq:prob1_eq5}
\dot{p}_w = c_v \frac{\partial^2 p_w}{\partial z^2},
\end{equation}
which is the basic differential equation of Terzaghi's consolidation theory \citep{terzaghi1943theoretical}. 
In this 1D setting, the boundary conditions can be specified as,
\begin{equation}
\label{eq:prob1_BC}
p_w = \hat{p}_w \text{ on } \partial \mathcal{B}_p
\: \: ; \: \:
-\frac{k}{\mu_w}\frac{\partial p_w}{\partial z} = \hat{w} \text{ on } \partial \mathcal{B}_w, 
\end{equation}
where $\hat{p}_w$ and $\hat{w}$ are the prescribed pore water pressure and Darcy's velocity on the Dirichlet ($\partial \mathcal{B}_p$) and Neumann ($\partial \mathcal{B}_w$) boundaries, respectively, that satisfy $\partial \mathcal{B} = \overline{\partial \mathcal{B}_p \cup \partial \mathcal{B}_w}$ and $\emptyset = \partial \mathcal{B}_p \cap \partial \mathcal{B}_w$. 
For model closure, the initial condition can be imposed as,
\begin{equation}
\label{eq:prob1_IC}
p_{w} = p_{w0},
\end{equation}
at time $t = 0$. 
Then, by following the standard weighted-residual procedure, the weak form of Eq.~\eqref{eq:prob1_eq5} can be obtained by multiplying it with a test function $\xi$ and integrating over $\mathcal{B}$:
\begin{equation}
\label{eq:prob1_weak}
\int_{\mathcal{B}} \xi \dot{p}_w \: \mathrm{d}z
+
\int_{\mathcal{B}} \frac{\partial \xi}{\partial z}c_v\frac{\partial p_w}{\partial z} \: \mathrm{d}z
+
\left. \xi K \hat{w} \right|_{\partial \mathcal{B}_w}
= 0.
\end{equation}

\subsection{Model Problem (2)}
\label{app:model_prob2}
Model Problem (2) focuses on a body of water-saturated soil $\mathcal{B} \subset \mathbb{R}^3$ where the solid and fluid phases are intrinsically compressible. 
This leads to the simplified effective stress equation: $\tensor{\sigma} = \tensor{\sigma}' - B p_w \tensor{I}$, such that balance of linear momentum [i.e., Eq.~\eqref{eq:bal_mom}] for a given saturated soil can be written as,
\begin{equation}
\label{eq:prob2_eq1}
\diver{(\tensor{\sigma}' - B p_w \tensor{I})} + \rho \vec{g} = \vec{0},
\end{equation}
while resulting in the following relation from Eqs.~\eqref{eq:mean_p} and \eqref{eq:lin_elas}: $(1 - \phi) p_s = -K \tr{(\tensor{\varepsilon})} + (B - \phi) p_w$.
By following \citet{sun2013stabilized}, if we assume that the change of porosity at infinitesimal time is small compare to those of intrinsic pressures, its total time derivative reads,
\begin{equation}
\label{eq:prob2_eq2}
(1-\phi) \dot{p}_s = -K \diver{\dot{\vec{u}}} + (B - \phi) \dot{p}_w.
\end{equation}
Since Eq.~\eqref{eq:bal_mass} can be simplified as, 
\begin{equation}
\label{eq:prob2_eq3}
\frac{1-\phi}{K_s} \dot{p}_s + \frac{\phi}{K_w}\dot{p}_w + \diver{\dot{\vec{u}}} + \diver{\vec{w}} = 0,
\end{equation}
substituting Eq.~\eqref{eq:prob2_eq2} into Eq.~\eqref{eq:prob2_eq3} gives the following expression for the mass balance of a saturated soil considered herein:
\begin{equation}
\label{eq:prob2_eq4}
\frac{1}{M}\dot{p}_w + B \diver{\dot{\vec{u}}} + \diver{\vec{w}} = 0,
\end{equation}
which is identical to the form that can be found in \citep{coussy2004poromechanics,belotserkovets2011thermoporoelastic,siriaksorn2018u}. 
Notice that by taking the displacement $\vec{u}$ and pore water pressure $p_w$ as the prime variables, Eqs.~\eqref{eq:prob2_eq1} and \eqref{eq:prob2_eq4} together form a so-called $\vec{u}$-$p_w$ formulation that requires mixed finite element methods to solve this model problem numerically. 
This requires considering the traction $\vec{t}$ and Darcy's velocity $\vec{w}$ as the secondary variables, such that the Dirichlet (displacement boundary $\partial \mathcal{B}_u$ and pore water pressure boundary $\partial \mathcal{B}_p$) and Neumann (traction boundary $\partial \mathcal{B}_t$ and water flux boundary $\partial \mathcal{B}_w$) boundary conditions for the governing partial differential equations can be specified as, 
\begin{equation}
\label{eq:prob2_BC}
\begin{dcases}
\vec{u} = \hat{\vec{u}} &\text{on } \partial \mathcal{B}_u \\
p_w = \hat{p}_w &\text{on } \partial \mathcal{B}_p
\end{dcases}
\: \: ; \: \:
\begin{dcases}
\tensor{\sigma} \cdot \vec{n} = \hat{\vec{t}} &\text{on } \partial \mathcal{B}_t \\
-\vec{w} \cdot \vec{n} = \hat{w} &\text{on } \partial \mathcal{B}_w
\end{dcases},
\end{equation}
where $\vec{n}$ indicates the outward unit normal on the boundary surface $\partial \mathcal{B}$ that satisfies:
\begin{equation}
\label{eq:prob2_BC_specs}
\partial \mathcal{B} = \overline{\partial \mathcal{B}_u \cup \partial \mathcal{B}_t} = \overline{\partial \mathcal{B}_p \cup \partial \mathcal{B}_w}
\: \: ; \: \:
\emptyset = \partial \mathcal{B}_u \cap \partial \mathcal{B}_t = \partial \mathcal{B}_p \cap \partial \mathcal{B}_w.
\end{equation}
To close the problem, the initial conditions are imposed as:
\begin{equation}
\label{eq:prob2_IC}
\vec{u} = \vec{u}_0 \: \: ; \: \: p_w = p_{w0},
\end{equation}
at time $t = 0$. 
Then, applying the standard weighted-residual procedure yields the weak form of Eqs.~\eqref{eq:prob2_eq1} and \eqref{eq:prob2_eq4} as follows:
\begin{equation}
\label{eq:prob2_weak}
\begin{dcases}
\int_{\mathcal{B}} \grad{\vec{\eta}} : \tensor{\sigma}' \: \mathrm{d}V - \int_{\mathcal{B}} B p_w \diver{\vec{\eta}} \: \mathrm{d} V - \int_{\mathcal{B}} \vec{\eta} \cdot \rho \vec{g} \: \mathrm{d}V - \int_{\partial \mathcal{B}_t} \vec{\eta} \cdot \hat{\vec{t}} \: \mathrm{d} \Gamma = 0, \\
\int_{\mathcal{B}} \xi \frac{1}{M} \dot{p}_w \: \mathrm{d}V + \int_{\mathcal{B}} \xi B \diver{\dot{\vec{u}}} \: \mathrm{d}V - \int_{\mathcal{B}} \grad{\xi} \cdot \vec{w} \: \mathrm{d}V - \int_{\partial \mathcal{B}_w} \xi \hat{w} \: \mathrm{d}\Gamma = 0,
\end{dcases}
\end{equation}
where $\vec{\eta}$ and $\xi$ are the test functions that correspond to the trial functions for the solution variables $\vec{u}$ and $p_w$, respectively.

\subsection{Model Problem (3)}
\label{app:model_prob3}
The third model problem considers a body of unsaturated soil $\mathcal{B} \subset \mathbb{R}^3$ with incompressible phase constituents. 
Even though unsaturated soil consists of three different phase constituents (e.g., solid, water, and air), the passive gas assumption outlined in Section \ref{sec:cont_rep} simplifies the problem using the pseudo-three-phase formulation, as it is not necessary to solve for the air phase explicitly. 
While this simplification may not fully capture all the complexities of the interactions among the three constituents, it makes computational process more efficient, while being reasonably accurate for near-surface applications \citep{white2016block, laloui2003solid, song2017strain}. 
In this case, note that the total stress can be expressed as: $\tensor{\sigma} = \tensor{\sigma}' - S^w p_w \tensor{I}$, which allows us to re-write Eq.~\eqref{eq:bal_mom} as the following:
\begin{equation}
\label{eq:prob3_eq1}
\diver{(\tensor{\sigma}' - S^w p_w \tensor{I})} + \rho \vec{g} = \vec{0},
\end{equation}
while using $K_s = K_w \approx \infty$, the mass balance equation in Eq.~\eqref{eq:bal_mass} can be re-written as,
\begin{equation}
\label{eq:prob3_eq2}
\phi \dot{S}^w + S^w \diver{\dot{\vec{u}}} + \diver{\vec{w}} = 0. 
\end{equation}
Similar to Appendix \ref{app:model_prob2}, by taking $\lbrace \vec{u}, p_w \rbrace$ as the prime variables and imposing the Neumann boundary conditions similar to Eqs.~\eqref{eq:prob2_BC}, the weak form of Eqs.~\eqref{eq:prob3_eq1} and \eqref{eq:prob3_eq2} can be written as,
\begin{equation}
\label{eq:prob3_weak}
\begin{dcases}
\int_{\mathcal{B}} \grad{\vec{\eta}} : \tensor{\sigma}' \: \mathrm{d}V - \int_{\mathcal{B}} S^w p_w \diver{\vec{\eta}} \: \mathrm{d} V - \int_{\mathcal{B}} \vec{\eta} \cdot \rho \vec{g} \: \mathrm{d}V - \int_{\partial \mathcal{B}_t} \vec{\eta} \cdot \hat{\vec{t}} \: \mathrm{d} \Gamma = 0, \\
\int_{\mathcal{B}} \xi \phi \dot{S}^w \: \mathrm{d}V + \int_{\mathcal{B}} \xi S^w \diver{\dot{\vec{u}}} \: \mathrm{d}V - \int_{\mathcal{B}} \grad{\xi} \cdot \vec{w} \: \mathrm{d}V - \int_{\partial \mathcal{B}_w} \xi \hat{w} \: \mathrm{d}\Gamma = 0.
\end{dcases}
\end{equation}

\end{appendices}

\end{document}